\documentclass[12pt]{article}
\usepackage{threeparttable}
\usepackage{amsmath}
\usepackage{enumerate}
\usepackage{natbib}
\usepackage{geometry}
\usepackage{url} 
\usepackage{graphicx,color,epstopdf}
\usepackage{amsfonts}
\usepackage{booktabs,multirow} 
\usepackage{subfig}
\usepackage{amsthm}
\usepackage{xspace}
\usepackage{url}

\usepackage{xcolor}
\usepackage[linesnumbered,ruled,vlined]{algorithm2e}
\SetKwInput{KwInput}{Input}                
\SetKwInput{KwOutput}{Output}              
\usepackage{setspace}
\definecolor{orange}{rgb}{1,0.5,0}
\usepackage[colorlinks,linkcolor=orange,anchorcolor=red,citecolor=orange]{hyperref}
\def\t{\top}

\newtheorem{theorem}{Theorem}

\newtheorem{lemma}{Lemma}
\newtheorem{remark}{Remark}

\newtheorem{assumption}{Assumption}

\begin{document}
	\date{}
		\title{Byzantine-robust distributed one-step estimation}
		\begin{center}
			{\LARGE\bf Byzantine-robust distributed one-step estimation}
			\\ \vskip0.6cm
			Chuhan Wang$^{1}$, Xuehu Zhu$^{2}$ and Lixing Zhu$^{1,3}\footnote{the research described herewith was supported by a grant from the Natural Scientific Foundation of China. The first two authors are the co-first authors.}$
			\\
		\textit{ $^1$ Center for Statistics and Data Science, Beijing Normal University, Zhuhai, China\\
				$^2$ School of Mathematics and Statistics, Xi’an Jiaotong University, Xi’an, China\\
				$^3$ Department of Mathematics,  Hong Kong Baptist University, Hong Kong, China }
		\end{center}

	\begin{abstract}
		 This paper proposes a Robust One-Step  Estimator (ROSE) to solve the Byzantine failure problem in distributed M-estimation when a moderate fraction of node machines experience Byzantine failures.
		To define   ROSE, the algorithms use the robust Variance Reduced Median Of the Local (VRMOL) estimator  to determine the initial parameter value for iteration,  and communicate between the node machines and the central processor in the Newton-Raphson iteration procedure to derive the robust VRMOL estimator of the gradient, and the Hessian matrix so as to obtain the final estimator.   ROSE has higher asymptotic relative efficiency than general median estimators without increasing the order of computational complexity. Moreover, this estimator can also cope with the problems involving anomalous or missing samples on the central processor. We prove the asymptotic normality when the parameter dimension $p$ diverges as the sample size goes to infinity, and  under weaker assumptions,  derive the  convergence rate.
		Numerical simulations and a real data application are conducted to evidence the effectiveness and robustness of  ROSE.
	\end{abstract}
\noindent%
{\it Keywords:} Distributed inference; Byzantine-robustness; variance reduction; M-estimaton.
	\section{Introduction}\label{sec1}
The growing scale of datasets due to technological advancements and decreasing data collection costs has made storing and processing massive amounts of data on a single computer increasingly challenging. Moreover, concerns over data privacy may make it inappropriate for certain machines to transmit data directly to others, necessitating the use of data processing techniques that avoid exposing sensitive information. Distributed statistical inference is a widely accepted solution to such challenges, where multiple machines are used to store and analyze data through the division of large tasks into smaller ones processed in parallel across different nodes. This approach enhances processing speed and system availability while reducing computing costs by utilizing inexpensive computer resources to complete large-scale computing tasks. 

In general, a distributed computing framework involves storing data on multiple machines with one central processor maintaining and updating target parameters, whereas other node machines perform local calculations based on their own data before transmitting essential information to the central processor. This framework faces two significant challenges. The first  is how to find a balance between computational accuracy and communication costs. The second is posed by the vulnerability of working machines and communication channels. A machine may experience a Byzantine failure, as defined by \cite{1982the}, wherein faulty machines compute some statistics into arbitrary or even contradictory values before transferring them to the central processor.  Byzantine failures may occur due to either machine failures or strong data heterogeneity  on some machines. 
For the former, some one-shot methods require only one communication between the central processor and the node machines, but perform poorly with nonlinear models.  See, e.g.,  \cite{2013Communication}; \cite{Jonathan2016On}; \cite{2017Communication}; \cite{2017Computational}; \cite{2018Distributed} and \cite{2020Learning}.  Iterative algorithms require multiple rounds of communication, leading to slow computation such as \cite{2017Efficient}; \cite{Jordan2019Communication} and \cite{2021Communication}. 
The latter requires robust estimation methods to  address the influence of Byzantine machines on estimators under a centralized framework.  See e.g.,  \cite{2014Distributed}; \cite{2018Byzantine,2019Defending} and \cite{su2019securing}. The key is to develop distributed estimators whose performances are as close to the global solution as possible, as if the data were aggregated in the central processor while ensuring the robustness of the estimators.

This paper focuses on developing  a  robust one-step M-estimation algorithm. This algorithm only requires one round of iteration.  The target is to estimate 
the solution  $\boldsymbol{\theta^*}$  of Equation (\ref{Me1}):
\begin{align}\label{Me1}
	\boldsymbol{\theta^*}=\underset{\boldsymbol{\theta} \in\boldsymbol{\Theta} }{\operatorname{argmin}} \mathbb{E}_{\boldsymbol{X}}[f(\boldsymbol{X}, \boldsymbol{\theta})]=\underset{\boldsymbol{\theta} \in\boldsymbol{\Theta} }{\operatorname{argmin}} \int_{\mathcal{X}}f(\boldsymbol{X}, \boldsymbol{\theta})dF(\boldsymbol{x}),
\end{align}
where  $\boldsymbol{\theta}=(\theta_1,\theta_2,\cdots,\theta_p)^{\top}$ is a $p$-dimensional vector in the parameter space $\boldsymbol{\Theta}\in \mathbb{R}^p$, and $\boldsymbol{X}=(X_1,X_2,\cdots,X_q)^{\top}$ is a $q$-dimensional random vector drawn from the sample space
$\mathcal{X}\subset\mathbb{R}^q$. Here, $F(\boldsymbol{x})$ represents the distribution function of  $\boldsymbol{X}$, $f$ is the general convex loss function, and $\mathbb{E}_{\boldsymbol{X}}(\cdot)$ stands for the expectation over $\boldsymbol{X}$.  We consider a distributed setting where $N$ independent and identically distributed (i.i.d.) samples, $\{\boldsymbol{X}_1,\boldsymbol{X}_2,\cdots,\boldsymbol{X}_N\}$, are  evenly distributed over $m$ machines  $\{\mathcal{I}_1,\mathcal{I}_2,\cdots,\mathcal{I}_m\}$, with each machine having $n$ samples. Under a Byzantine distributed framework, a small proportion of machines, indexed by a set $\mathcal{B}\subseteq\{1,2,3,\cdots,m\}$, where $\alpha_n m$ elements belong to $\mathcal{B}$ and $\alpha$ denotes the proportion of Byzantine machines. When these Byzantine machines compute statistics using their respective data, they return $*$ which denotes an arbitrary value, whereas normal machines output the correct statistics.

Note that when all machines are normal, existing efficient parameter estimation methods often use gradient descent to update parameters, such as the ``Communication-efficient Surrogate Likelihood (CSL)" framework proposed by \cite{Jordan2019Communication} that computes the average of the gradients from the node machines and updates the parameter estimator through iterative improvement.  As an alternative, \cite{huo2019distributed} propose a one-step approach that the central processor employs the averages  of these gradients and Hessians from the node machines to adjust the initial estimator. This one-step approach only requires a single extra round of communication following the calculation of the initial estimator. Under the presence of Byzantine machines, the  ``median-of-means" (MOM) estimator has been widely studied in the literature such as \cite{2018Byzantine}; \cite{minsker2015geometric,minsker2019distributed}; \cite{2019Regularization} and \cite{Guillaume2020Robust}. It computes the local sample mean for each node machine and uses the central processor to calculate the median. The MOM estimator is robust, even when the proportion of Byzantine machines is nearly 50 percents, but its asymptotic efficiency relative to the mean estimator is only $2/\pi\approx 0.637$. 
\cite{tu2021variance} apply the Composite Quantile Regression (CQR) estimator proposed by \cite{2008Composite} to the gradient in a distributed framework when the initial estimator and Hessian matrix are computed by a certain normal machine using its own data as the CQR estimator can reach an asymptotic efficiency reaching $3/\pi\approx 0.955$ for linear regression models. This method relies on a precise initial parameter estimator and adequate information about the Hessian matrix to obtain an accurate estimator. 
In the conducted numerical studies, we found that the differences between the gradient estimator and the actual gradient value significantly affects its estimation stability.  Further,  as this method utilizes the samples on the central processor, when  the central processor is a Byzantine machine with anomalous samples on it or with no samples,  their estimator may not work well or not be implementable.

\subsection{Our contributions}
The primary contributions of this paper are as follows. First, we apply the concept of the composite quantile regression (CQR) estimation 
to construct a more general estimator called the Varianced Reduced Median Of Local (VRMOL) estimator that aggregates estimators on node machines for the M-estimator. We argue that when
$n$, $m$ and $p$ satisfy certain conditions, the VRMOL estimator employed as the initial parameter value can achieve higher asymptotic efficiency relative to calculating each node machine estimator's median directly. Second, we use this robust VRMOL estimator to estimate the gradient and Hessian matrix such that we can use a similar method to \cite{huo2019distributed} to correct our initial estimator when only one iteration is used. 
We show that this approach works in diverging
$p$ scenarios under some circumstances. 
Worth noting, since we only perform one iteration, non-convergence caused by gradients not approaching zero (e.g., \cite{Jordan2019Communication}, \cite{tu2021variance}) is not a concern. The numerical studies also support this claim. Third, we propose an algorithm that can produce consistent parameter estimators even when samples on the central processor are anomalous or unavailable. 
Finally, we derive the asymptotic normality of the Robust One-Step  Estimator (ROSE) under the rate conditions $\alpha_n=o(1/\sqrt{pm}), pm\log n/n=o(1)$, $p\log^2 n/m=o(1)$ and $p^2\sqrt{m}\log^{3/2} n/n=o(1)$ so that  the estimator can be used for statistical inference. Furthermore, we derive the estimator's convergence rate under weaker rate assumptions than the above.

The rest of the paper is organized as follows. Section \ref{sec2} describes the motivation and the properties of the VRMOL estimator. 
Section \ref{sec3} introduces the ROSE algorithm   so that  the VRMOL estimator can be applied to estimating the initial parameter value, the gradient, and the Hessian matrix simultaneously. Section \ref{sec4} discusses the asymptotic normality and the convergence rate of  ROSE.  Section \ref{sec5} includes numerical results on synthetic and real data,  and Section \ref{sec6} presents the concluding remarks and further discussions. Section \ref{sec7} contains all the assumptions. The technical proofs of the results are postponed to  Supplementary Material.

\subsection{Notations}
Give some notations first. For any vector $\boldsymbol{v}=(v_1,v_2,\cdots,v_p)^{\top}$, define $\|\boldsymbol{v}\|=(\sum_{l=1}^{p}v_l^2)^{1/2}$. Let 
$\boldsymbol{v}^{\otimes 2}=\boldsymbol{v}\boldsymbol{v}^{\top}$. Define $B(\boldsymbol{v},r)=\{\boldsymbol{w}\in\mathbb{R}^p:\|\boldsymbol{w}-\boldsymbol{v}\|\leq r\}$. For two vectors $\boldsymbol{u}\in\mathbb{R}^p$ and $\boldsymbol{v}\in \mathbb{R}^p$, let $\langle\boldsymbol{u},\boldsymbol{v}\rangle=\sum_{l=1}^{p}u_lv_l$. For a matrix $\mathbf{M}$,  denote  $\|\mathbf{M}\|=\sup_{\|\boldsymbol{a}\|=1}{\|\boldsymbol{a}^{\top}\mathbf{M}\|}$ as the spectral norm, which is equal to the largest eigenvalue of $\mathbf{M}$ if $\mathbf{M}$ is a square matrix.  $\lambda_{\max}(\mathbf{M})$ and $\lambda_{\min}(\mathbf{M})$ are the largest  and  smallest eigenvalues of $\mathbf{M}$ respectively. Let $\operatorname{diag}(\mathbf{M})$ be the vector which is composed of entries on the diagonal of $\mathbf{M}$. Denote $\operatorname{Tr}(\mathbf{M})$ as the trace of $\mathbf{M}$. Define $\mathbf{N}(0,1)$ as the standard normal distribution, $\Psi(x)=\mathbb{P}(\mathbf{N}(0,1)\leq x)$ and $\psi(x)=\frac{1}{\sqrt{2\pi}}\exp(-\frac{x^2}{2})$, which are the distribution and density functions of the standard normal distribution respectively. Let $F^{-1}(y)$ be the smallest $x\in\mathbb{R}$ satisfing $F(x)\geq y$. Denote $\nabla$ and $\nabla^2$ as the gradient and Hessian operators. For $f(\boldsymbol{\theta})=f(\theta_1,\theta_2,\cdots,\theta_p)$, $\nabla_{\theta_l}f$ stands for  the partial derivatives with respect to $\theta_l$. Denote $\stackrel{\mathbb{P}}{\to}$ and $\stackrel{d}{\to}$ as convergence in probability and  in distribution respectively. Define $a_n\asymp b_n$ if $a_n=O(b_n)$ and $b_n=O(a_n)$.  Let $\mathbb{I}(\cdot)$ be the indicator function. Let $\mathcal{A}^c$ be the complement of event $\mathcal{A}$. For any positive integer $N$, denote $[N]$ as the index set $\{1,2,\cdots,N\}$. For $X_1,X_2,\cdots,X_m$, let $\operatorname{med}\{X_j,j\in[m]\}$ be the median of $\{X_1,X_2,\cdots,X_m\}$. Write the $m$ machines as $\mathcal{I}_1$, $\mathcal{I}_2$, $\cdots$, $\mathcal{I}_m$. Here $\mathcal{I}_1$ will act as the central processor, and the others  are node machines.

	\section{The VRMOL estimator}\label{sec2}
We subsequently investigate its asymptotic properties in the following subsections.
	\subsection{Motivation}
	Consider a general M-estimation problem in a distributed framework, where the objective is to minimize the loss function $f(\boldsymbol{X},\boldsymbol{\theta})$, in expectation, with the target parameter $\boldsymbol{\theta^*}$ defined in Equation (\ref{Me1}). The entire dataset has $N=mn$ i.i.d. observations $\{\boldsymbol{X}_i\}$, $i\in [N]$, and is evenly distributed on $m$ machines $\{\mathcal{I}_j\}_{j=1}^m$, with each machine containing $n$ observations. We use $i\in \mathcal{I}_j$ to represent that $\boldsymbol{X}_i$ belongs to machine $\mathcal{I}_j$.
	For any $j\in[m]$,  define
	\begin{align}\label{th1}
		F_j(\boldsymbol{\theta})=\frac{1}{n}\sum_{i\in \mathcal{I}_j}f(\boldsymbol{X}_i,\boldsymbol{\theta})
		\quad \text{and} \quad \boldsymbol{\hat{\theta}}_j=(\hat{\theta}_{j1},\cdots,\hat{\theta}_{jp})=\underset{\boldsymbol{\theta}\in\boldsymbol{\Theta}}{\operatorname{argmin}}F_j(\boldsymbol{\theta}).
	\end{align}
	Therefore, for the global loss function, we have
	\begin{align}\label{th2}
		F(\boldsymbol{\theta})=\frac{1}{m}\sum_{j=1}^{m}F_j(\boldsymbol{\theta})=\frac{1}{N}\sum_{i=1}^{N}f(\boldsymbol{X}_i,\boldsymbol{\theta}).
	\end{align}
	
	Additionally, for convenience, we define $F_{\mu}(\boldsymbol{\theta})=\mathbb{E}\{f(\boldsymbol{X},\boldsymbol{\theta})\}$, where $\boldsymbol{X}$ shares the same distribution with $\boldsymbol{X}_i$. Then, the target parameter is given by $\boldsymbol{\hat{\theta}^*}={\operatorname{argmin}}_{\boldsymbol{\theta}\in\boldsymbol{\Theta}}{F(\boldsymbol{\theta})}$. In a distributed setting, directly obtaining the value of $\boldsymbol{\hat{\theta}^*}={\operatorname{argmin}}_{\boldsymbol{\theta}\in\boldsymbol{\Theta}}{F(\boldsymbol{\theta})}$ is not possible. Therefore, \cite{2013Communication} propose using the average of $\boldsymbol{\hat{\theta}}_j$ for $j\in[m]$ to estimate $\boldsymbol{\theta^*}$. The corresponding estimator achieves a mean squared error of order $O(1/mn)$ under certain assumptions. However, in a Byzantine distributed setting, we must incorporate some robust methods to estimate the parameter. An idea is to use the median of $\boldsymbol{\hat{\theta}}_j$, i.e., $\boldsymbol{\hat{\theta}}_{med}=(\hat{\theta}_{med,1},\cdots,\hat{\theta}_{med,p})^{\t}$, where $\hat{\theta}_{med,l}=\operatorname{med}\{\hat{\theta}_{jl},j\in[m]\}$, $l\in[p]$, instead of the average to estimate $\boldsymbol{\hat{\theta}^*}$. However, although this estimator is easy to implement, its asymptotic relative efficiency may be very low. For example,  it is not difficult to prove that for any $l\in[p]$,  the mean estimator $\boldsymbol{\hat{\theta}}_j=\frac{1}{n}\sum_{i\in\mathcal{I}_j}\boldsymbol{X}_i$ has the normal weak limit:
	\begin{align*}
		\frac{\sqrt{mn}(\hat{\theta}_{med,l}-\mathbb{E}(X_l))}{\sigma(X_l)}\stackrel{d}{\to}\mathbf{N}(0,\frac{\pi}{2}),
	\end{align*}
	where $\sigma(X_l)$ is the  standard deviation of $X_l$. Then the asymptotic efficiency relative to the mean estimator is ${2}/{\pi}\approx 0.637$, which is far away from $1$. To improve the estimation effectiveness, we now  introduce a new estimator. 

	Consider a quantile loss function defined as $\mathcal{L}_{\kappa}(z)=z(\kappa-\mathbb{I}(z\leq 0))$ with $0<\kappa<1$. Assume that $Y_j,j\in[m]$  are independent and identically distributed random variables with a symmetric distribution. Let  $\mathbb{E}(Y_j)=\mu_Y={\operatorname{argmin}}_{y\in \mathbb{R}}\mathbb{E}[\mathcal{L}_{1/2}(Y-y)]$, and $G_{Y}(y)=\mathbb{P}(Y\leq y)$ be the distribution function. The sample median, $\hat{Y}_{med}=\operatorname{med}\{Y_j,j\in[m]\}$, can be
	a consistent estimator of ${\operatorname{argmin}}_{y\in \mathbb{R}}\mathbb{E}[\mathcal{L}_{1/2}(Y-y)]$, where $Y$ has the same distribution as $Y_j$. It is worth noting that not only the median of $\operatorname{med}\{Y_j,j\in[m]\}$ is an unbiased estimator of $\mu_Y$, but also the average of the $\kappa$ quantile and $1-\kappa$ quantile of $\{Y_j,j\in[m]\}$ provides an unbiased estimator of $\mu_Y$. This is because $\mu_Y=[G_Y^{-1}(\kappa)+G_Y^{-1}(1-\kappa)]/2$. Consequently, an estimator of  $\mu_Y$
	can be obtained by averaging the estimators of $G_Y^{-1}(\kappa)={\operatorname{argmin}}_{y\in \mathbb{R}}\mathbb{E}[\mathcal{L}_{\kappa}(Y-y)]$ and $G_Y^{-1}(1-\kappa)={\operatorname{argmin}}_{y\in \mathbb{R}}\mathbb{E}[\mathcal{L}_{1-\kappa}(Y-y)]$. This estimator proves reliable when the number of outliers is fewer than $\min\{m\kappa,m(1-\kappa)\}$. The method of averaging quantiles can yield a smaller variance compared to the median estimator. However, when numerous pairs of quantile estimators are used, computing the average may become complex. For symmetric distributions of $Y$, $\mu_Y$ can also be expressed as $G_Y^{-1}(\kappa)-[G_Y^{-1}(\kappa)-G_Y^{-1}(1/2)]$, which is equivalent to ${\operatorname{argmin}}_{y\in \mathbb{R}}\mathbb{E}[\mathcal{L}_{\kappa}(Y-[G_Y^{-1}(\kappa)-G_Y^{-1}(1/2)]-y)]$. Thus, for $\kappa_1,\kappa_2,\cdots,\kappa_K$, $\mu_Y={\operatorname{argmin}}_{y\in \mathbb{R}}\mathbb{E}[\sum_{k=1}^{K}\mathcal{L}_{\kappa_k}(Y-[G_Y^{-1}(\kappa_k)-G_Y^{-1}(1/2)]-y)]$. To employ several quantiles in reducing the estimator's variance, we can directly contemplate minimizing the following expectation:
\begin{align*}
	\mathcal{E}(y):=\mathbb{E}\left[\sum_{k=1}^{K}\mathcal{L}_{\kappa_k}(Y-[G_{Y}^{-1}(\kappa_k)-G_{Y}^{-1}(1/2)]-y)\right].
\end{align*}
Now we discuss how to estimate ${\operatorname{argmin}}_{y\in \mathbb{R}}\mathcal{E}(y)$. Although it is possible to directly obtain the estimator by minimizing the sample form of $\mathcal{E}(y)$, the computational cost is considerable due to the loss function's non-smoothness. To alleviate the computational burden, a Newton-Raphson method is proposed under the assumption that the density function $g_Y(y)$ of $Y$ is continuous. Using an easily estimable and approximate initial estimator for $\mu_Y$, such as $\hat{Y}_{med}$, the Newton-Raphson iteration can help estimate $\mu_Y$.
Here we need to take advantage of the gradient
\begin{align*}
	\frac{\partial \mathcal{E}(y)}{\partial y}=\sum_{k=1}^{K}\mathbb{E}\left[\mathbb{I}(Y\leq y+[G_{Y}^{-1}(\kappa_k)-G_{Y}^{-1}(1/2)])-\kappa_k\right]
\end{align*}
and the second derivative
\begin{align*}
	\frac{\partial^2 \mathcal{E}(y)}{\partial y^2}=\sum_{k=1}^{K}g_Y(y+[G_{Y}^{-1}(\kappa_k)-G_{Y}^{-1}(1/2)]).
\end{align*}
Since we have yet to learn the distribution function $G_Y(y)$, it is still difficult to estimate the gradient and the second derivative. If $Y_j,j\in[m]$ are the estimators computed by different machines with a tractable limiting distribution, we can use the limiting distribution function $\hat{G}_Y$ instead of $G_Y$ to estimate $\mu_Y$.

 We present the following variance-reduced estimator using $\hat{Y}_{med}$ as the initial estimator:
\begin{align}\label{0301}
	\hat{\mu}_{Y}=\hat{Y}_{med}-\frac{\frac{1}{m}\sum_{k=1}^{K}\sum_{j=1}^{m}\left[\mathbb{I}(Y_j\leq \hat{Y}_{med}+[\hat{G}_{Y}^{-1}(\kappa_k)-\hat{G}_{Y}^{-1}(1/2)])-\kappa_k\right]}{\sum_{k=1}^{K}\hat{g}_Y(\hat{Y}_{med}+[\hat{G}_{Y}^{-1}(\kappa_k)-\hat{G}_{Y}^{-1}(1/2)])},
\end{align}
where $\hat{g}_Y$ is the density function corresponding to $\hat{G}_Y$. This is the crucial idea of reducing the variance of the median estimator, see e.g., \cite{2008Composite}. 


\begin{remark}
In this study, we employ $\hat{Y}_{med}$ as the initial estimator in (\ref{0301}). Alternative estimators, such as the coordinate-wise trimmed mean, may also be utilized.  If the convergence rate of an initial estimator  is the same as $\hat{Y}_{med}$,  the final estimator $\hat{\mu}_{Y}$  also has the same convergence rate as the estimator using $\hat{Y}_{med}$  as the initial one.
\end{remark}

\subsection{Application to M-estimation}
We adopt the variance reduction idea to the general M-estimation in a distributed framework to define an estimator $\boldsymbol{\hat{\theta}}_{vr}$.

Define
\begin{align*}
	\boldsymbol{\Sigma}(\boldsymbol{\theta})=\{\nabla^2F_{\mu}(\boldsymbol{\theta})\}^{-1}\mathbb{E}\{\nabla f(\boldsymbol{X},\boldsymbol{\theta})^{\otimes 2}\} \{\nabla^2F_{\mu}(\boldsymbol{\theta})\}^{-1},
\end{align*}
and let $\sigma_{l_1,l_2}(\boldsymbol{\theta})$ be the $(l_1,l_2)$-entry of $\boldsymbol{\Sigma}(\boldsymbol{\theta})$. To simplify the notations, define $\sigma_{l}^2(\boldsymbol{\theta})$ as the $(l,l)$-entry of $\boldsymbol{\Sigma}(\boldsymbol{\theta})$.
Recall $\mathcal{B} \subseteq\{1,2,3, \ldots, m\}$ the  subset of the indices and $\alpha_nm$ with $\alpha_n<1/2$ the number of elements in $\mathcal{B}$. Define
$$
\boldsymbol{\hat{\theta}}_{j}=\left\{\begin{array}{cl}
	\underset{\boldsymbol{\theta}\in\boldsymbol{\Theta}}{\operatorname{argmin}}F_j(\boldsymbol{\theta}) & j \notin \mathcal{B}, \\
	* & j \in \mathcal{B},
\end{array}
\right.
$$
where $*$ can be any value given by the Byzantine machine.	

To take the advantage of the algorithm provided by (\ref{0301}),  we utilize the asymptotic normality of the M-estimator. Recall that $\boldsymbol{\hat{\theta}}_j$ is the local M-estimator on the machine $\mathcal{I}_j$. Under Assumptions \ref{a1}-\ref{a4} and \ref{a8}, if $p^2\log^2 n/n=o(1)$, for each entry in $\boldsymbol{\hat{\theta}}_j=(\hat{\theta}_{j1},\cdots,\hat{\theta}_{jp})$, we can prove that for any $l\in[p]$ (see Lemma 15 in  Supplementary Material for details),
\begin{align*}
	\frac{\sqrt{n}(\hat{\theta}_{jl}-\theta_l^*)}{\sigma_l(\boldsymbol{\theta^*})}\stackrel{d}{\to}\mathbf{N}(0,1).
\end{align*}
However, since $\sigma_l(\boldsymbol{\theta^*}), l\in
[p]$ are unknown, we must estimate them before using the algorithm in (\ref{0301}). Here, we provide two methods for this purpose.
One  method only utilizes data from the central processor and ${\boldsymbol{\hat{\theta}}}_{med}:$
\begin{align}\label{sig2}
	&(\hat{\sigma}_1^2({\boldsymbol{\hat{\theta}}}_{med}),\cdots,\hat{\sigma}_p^2({\boldsymbol{\hat{\theta}}}_{med}))\nonumber
	\\&\ \ \ =\operatorname{diag}(\{\nabla^2F_1(\boldsymbol{\boldsymbol{\hat{\theta}}}_{med})\}^{-1}\frac{1}{n}\sum_{i\in\mathcal{I}_1}\{\nabla f(\boldsymbol{X}_i,\boldsymbol{\boldsymbol{\hat{\theta}}}_{med})-\nabla F_1(\boldsymbol{\hat{\theta}}_{med})\}^{\otimes 2} \{\nabla^2F_1(\boldsymbol{\boldsymbol{\hat{\theta}}}_{med})\}^{-1}),
\end{align}
where $\boldsymbol{\hat{\theta}}_{med}=(\hat{\theta}_{med,1},\cdots,\hat{\theta}_{med,p})^{\t}$ and $\hat{\theta}_{med,l}=\operatorname{med}\{\hat{\theta}_{jl},j\in[m]\}$, $l\in[p]$. However, if the data on the central processor is anomalous, this estimator may deviate significantly from the true value. To solve this problem,  another estimator, which requires the node machines to transmit  variance estimators $\hat{\sigma}_{jl}^2(\boldsymbol{\hat{\theta}}_j), l\in[p]$ while transmitting $\boldsymbol{\hat{\theta}}_j$ to the central processor, is provided as:
\begin{align}\label{sig4}
	\hat{\sigma}{'}_l^2=\operatorname{med}\{\hat{\sigma}_{jl}^2({\boldsymbol{\hat{\theta}}}_j),j\in[m]\}, l\in[p],
\end{align}
where $(\hat{\sigma}_{j1}^2({\boldsymbol{\hat{\theta}}}_j),\cdots,\hat{\sigma}_{jp}^2({\boldsymbol{\hat{\theta}}}_j))=\operatorname{diag}(\{\nabla^2F_j({\boldsymbol{\hat{\theta}}}_j)\}^{-1}\frac{1}{n}\sum_{i\in\mathcal{I}_j}\{\nabla f(\boldsymbol{X}_i,{\boldsymbol{\hat{\theta}}}_j)-\nabla F_j({\boldsymbol{\hat{\theta}}}_j)\}^{\otimes 2} \{\nabla^2F_j({\boldsymbol{\hat{\theta}}}_j)\}^{-1})$.

Then we have the following lemma to state the consistency of these two estimators.
\begin{lemma}\label{la9} Suppose that 	$\log n/m=o(1)$, $\alpha_n=o(1)$ and $p^2\log^{2} n/n=o(1)$. Then under Assumptions \ref{a1}-
	\ref{a9} and \ref{a11}-\ref{a12},  for any $l\in[p]$,
\begin{itemize}
\item	
[(1)] If we have prior information that $\mathcal{B}\subseteq\{2,3,\cdots,m\}$, then
	\begin{align*}
		|\hat{\sigma}_l(\boldsymbol{\hat{\theta}}_{med})-\sigma_l(\boldsymbol{\theta^*})|=O_p(p/\sqrt{n}).
	\end{align*}

\item [(2)] {For $\mathcal{B}\subseteq\{1,2,3,\cdots,m\}$,
\begin{align*}
	|\hat{\sigma}_l{'}-\sigma_l(\boldsymbol{\theta^*})|=O_p(p/\sqrt{n}).
\end{align*}}
\end{itemize}
\end{lemma}
For convenience, in the following discussions of this section, we assume $\mathcal{B}\subseteq\{2,3,\cdots,m\}$ and use $\hat{\sigma}_l(\boldsymbol{\hat{\theta}}_{med})$ as the estimator of $\sigma_l(\boldsymbol{\theta^*})$, but these discussions are also available when $\mathcal{B}\subseteq\{1,2,3,\cdots,m\}$ if we use $\hat{\sigma}_l{'}$ as the estimator of $\sigma_l(\boldsymbol{\theta^*})$.
Under Assumption \ref{a9}, we can use $\hat{G}_l(t)=\Psi(\sqrt{n}(t-\theta^*_l)/\hat{\sigma}_l(\boldsymbol{\hat{\theta}}_{med}))$ and $\hat{g}_l(t)=\sqrt{n}\psi(\sqrt{n}(t-\theta^*_l)/\hat{\sigma}_l(\boldsymbol{\hat{\theta}}_{med}))/\hat{\sigma}_l(\boldsymbol{\hat{\theta}}_{med})$ to estimate the cumulative distribution   and  density functions of $\hat{\theta}_{jl}$.
Note that $\hat{G}_l^{-1}(\kappa_k)=\theta^*_l+\hat{\sigma}_l(\boldsymbol{\hat{\theta}}_{med})\Psi^{-1}(\kappa_k)/\sqrt{n}$. Then according to (\ref{0301}), we infer that
\begin{align}\label{03021}
	&\hat{\theta}_{med,l}-\frac{\frac{1}{m}\sum_{k=1}^{K}\sum_{j=1}^{m}\left[\mathbb{I}(\hat{\theta}_{jl}\leq \hat{\theta}_{med,l}+\hat{\sigma}_l(\boldsymbol{\hat{\theta}}_{med})\Psi^{-1}(\kappa_k)/\sqrt{n})-\kappa_k\right]}{\sum_{k=1}^{K}\sqrt{n}\psi(\sqrt{n}(\hat{\theta}_{med,l}-\theta^*_l)/\hat{\sigma}_l(\boldsymbol{\hat{\theta}}_{med})+\Psi^{-1}(\kappa_k))/\hat{\sigma}_l(\boldsymbol{\hat{\theta}}_{med})}\nonumber
	\\=&\hat{\theta}_{med,l}-\frac{\hat{\sigma}_l(\boldsymbol{\hat{\theta}}_{med})\sum_{k=1}^{K}\sum_{j=1}^{m}\left[\mathbb{I}(\hat{\theta}_{jl}\leq \hat{\theta}_{med,l}+\hat{\sigma}_l(\boldsymbol{\hat{\theta}}_{med})\Psi^{-1}(\kappa_k)/\sqrt{n})-\kappa_k\right]}{m\sqrt{n}\sum_{k=1}^{K}\psi(\Psi^{-1}(\kappa_k))}+o_p(1)\nonumber
	\\=: \, &\hat{\theta}_{vr,l}^K+o_p(1),
\end{align}
where the first equation is derived by using the fact that $\sqrt{n}(\hat{\theta}_{med,l}-\theta^*_l)=o_p(1)$. This result is stated in the following lemma.

\begin{lemma}\label{oMom}
	(Concentration of the median of the local M-estimators with Byzantine machines) Suppose Assumptions \ref{a1}-\ref{a9} and  \ref{a12} hold.
	There exists a constant $\tilde{C}_{\gamma}>0$ only depending on $\gamma$, where $\gamma$ can be any constant larger than $1$, such that
	$$
	\mathbb{P}\left\{|{\hat{\theta}}_{\operatorname{med},l}-\theta^*_l| \geq \frac{\gamma\tilde{C}_{\gamma}}{\sqrt{n}}\left(\alpha_{n}+\frac{p\log n}{\sqrt{n}}+\sqrt{\frac{\log n}{m}}\right)\right\}=o(n^{-\gamma}).
	$$
\end{lemma}
\begin{remark}
	Lemma \ref{oMom}  provides the convergence rate of the Median Of the Local(MOL) estimator. According to Theorem 1 in \cite{2013Communication}, when $p$ is fixed, if all machines are normal, the convergence rate of the estimator $\boldsymbol{\hat{\theta}}_{mean}$ obtained by directly averaging the local estimators of all machines is $O_p(1/\sqrt{mn}+1/n)$ that is only faster than the MOL estimator ${\boldsymbol{\hat{\theta}}_{\operatorname{med}}}=({\hat{\theta}}_{\operatorname{med},1},\cdots,{\hat{\theta}}_{\operatorname{med},p})$ with a rate of order ${\log n}$. {If we put all the normal data together and calculate a global estimator $\boldsymbol{\hat{\theta}^*}$, then when $p$ is divergent, the convergence rate of $\boldsymbol{\hat{\theta}^*}$  is $\sqrt{p/mn}$.  Compared to $\boldsymbol{\hat{\theta}^*}$, when $\alpha_n=o(\sqrt{\log n/m})$ and $p^2m\log n/n=o(1)$, the MOL estimator converges slightly slower with a factor  $\sqrt{\log n}$. }
\end{remark}
Now we provide the convergence rate of $\boldsymbol{\hat{\theta}}_{vr}=(\hat{\theta}_{vr,1}^K,\cdots,\hat{\theta}_{vr,p}^K)$, where $\hat{\theta}_{vr,l}^K, l\in[p]$ has been defined in (\ref{03021}). To make our conclusion clear,
let $\kappa_k=\frac{k}{K+1}$ and $\Delta_k=\Psi^{-1}(\kappa_k)=\Psi^{-1}(\frac{k}{K+1})$.
Then
\begin{align}\label{vr4}
	\hat{\theta}_{vr,l}^K=\hat{\theta}_{med,l}-\frac{\hat{\sigma}_l(\boldsymbol{\hat{\theta}}_{med})\sum_{k=1}^{K}\sum_{j=1}^{m}\left[\mathbb{I}(\hat{\theta}_{jl}\leq \hat{\theta}_{med,l}+\hat{\sigma}_l(\boldsymbol{\hat{\theta}}_{med})\Delta_k/\sqrt{n})-\frac{k}{K+1}\right]}{m\sqrt{n}\sum_{k=1}^{K}\psi(\Delta_k)}.
\end{align}

To investigate  the asymptotic properties of $\boldsymbol{\hat{\theta}}_{vr}$, we define a new matrix $\boldsymbol{\Sigma}_{vr}^{K}(\boldsymbol{\theta^*})$ in the following.
Let $(\xi_{l_1},\xi_{l_2})$ follow a zero mean bivariate normal distribution with $\operatorname{Var}(\xi_{l_1})=\operatorname{Var}(\xi_{l_2})=1$ and $\operatorname{Cov}(\xi_{l_1},\xi_{l_2})=\frac{\sigma_{l_1,l_2}(\boldsymbol{\theta^*})}{{\sigma_{l_1}(\boldsymbol{\theta^*})\sigma_{l_2}(\boldsymbol{\theta^*}})}$. Denote $\kappa_{k_1,k_2,K}^{l_1,l_2}=\mathbb{P}(\xi_{l_1}\leq\Delta_{k_1},\xi_{l_2}\leq\Delta_{k_2})$.
For $l_1\neq l_2$, define the $(l_1,l_2)$-entry of the matrix $\boldsymbol{\Sigma}_{vr}^{K}(\boldsymbol{\theta^*})$ as
$$\sigma_{vr,l_1,l_2}=\frac{\sum_{k_1=1}^K\sum_{k_2=1}^K(\kappa_{k_1,k_2,K}^{l_1,l_2}-\kappa_{k_1}\kappa_{k_2})}{\{\sum_{k=1}^{K}\psi(\Delta_{k})\}^2}\sigma_{l_1}(\boldsymbol{\theta^*})\sigma_{l_2}(\boldsymbol{\theta^*}),$$
where $\kappa_{k}=\frac{k}{K+1}$. For the $(l,l)$-entry $\sigma_{vr,l}^2$ in $\boldsymbol{\Sigma}_{vr}^{K}(\boldsymbol{\theta^*})$, noting that $\kappa_{k_1,k_2,K}^{l,l}=\mathbb{P}(\xi_{l}\leq\min\{\Delta_{k_1},\Delta_{k_2}\})=\min\{\kappa_{k_1},\kappa_{k_2}\}$,  then we have
\begin{align}\label{DK}
	\sigma_{vr,l}^2=\frac{\sum_{k_1=1}^K\sum_{k_2=1}^K(\min\{\kappa_{k_1},\kappa_{k_2}\}-\kappa_{k_1}\kappa_{k_2})}{\{\sum_{k=1}^{K}\psi(\Delta_{k})\}^2}\sigma_{l}^2(\boldsymbol{\theta^*})=:D_K\sigma_{l}^2(\boldsymbol{\theta^*}).
\end{align}

Similarly to (\ref{vr4}), by replacing $\hat{\sigma}_l(\boldsymbol{\hat{\theta}}_{med})$ with $\hat{\sigma}{'}_l$, we can define $\boldsymbol{\hat{\theta}}{'}_{vr}$. The following theorem states the convergence rate and asymptotic normality of $\boldsymbol{\hat{\theta}}_{vr}$ and $\boldsymbol{\hat{\theta}}{'}_{vr}$.

\begin{theorem}\label{key1}
	(Variance Reduced Median Of Local Estimator, VRMOL)  Suppose that $\mathcal{B}\subseteq\{2,3,\cdots,m\}$ and Assumptions \ref{a1}-\ref{a9} and \ref{a11}-\ref{a12} hold.  Replace $\alpha_nm$ estimators of $\boldsymbol{\hat{\theta}}_1$, $\cdots$,$\boldsymbol{\hat{\theta}}_m$ defined in (\ref{th1}) with arbitrary values.  $\boldsymbol{\hat{\theta}}_{vr}=(\hat{\theta}_{vr,1}^K,\cdots,\hat{\theta}_{vr,p}^K)^{\top}$ is defined in (\ref{03021}). When $\alpha_n=O(1/\log n)$, ${\log^3 n}/{m}=o(1)$ and ${p^2\log^2 n}/{n}=o(1)$,
	\begin{align*}
		\hat{\theta}_{vr,l}^K-	\theta_l^*=O_p\left(\frac{\alpha_n}{\sqrt{n}}+\sqrt{\frac{1}{mn}}+\frac{p\log n}{n}\right),
	\end{align*}
	and
	\begin{align*}
		\|\boldsymbol{\hat{\theta}}_{vr}-	\boldsymbol{\theta^*}\|=O_p\left(\frac{\alpha_n\sqrt{p}}{\sqrt{n}}+\sqrt{\frac{p}{mn}}+\frac{p^{3/2}\log n}{n}\right).
	\end{align*}
	When $p^2m\log^2 n/n=o(1)$,
	\begin{align}\label{t15}
		\frac{\sqrt{mn}(\hat{\theta}_{vr,l}^K-\theta^*_l)}{\sigma_{vr,l}}\stackrel{d}{\to}\mathbf{N}(0,1),
	\end{align}
	where
	$\sigma_{vr,l}^2=D_K\sigma_{l}^2(\boldsymbol{\theta^*})$ and $D_K$ is defined in (\ref{DK}).
	
		Moreover, when $\alpha_n=o(\max\{{1}/{\sqrt{mp}},{1}/(p\log n)\})$, $p^2\log^3 n/m=o(1)$ and ${p^3m\log^2 n}/n=o(1)$, for any constant vector $\boldsymbol{v}$ satisfying $\|\boldsymbol{v}\|=1$,
	\begin{align}\label{t14}
		\frac{\sqrt{mn}}{\sqrt{\boldsymbol{v}^{\top}\mathbf{\Sigma}_{vr}^{K}(\boldsymbol{\theta^*})\boldsymbol{v}}}\langle\boldsymbol{\hat{\theta}}_{vr}-	\boldsymbol{\theta^*},\boldsymbol{v}\rangle\stackrel{d}{\to}\mathbf{N}(0,1).
	\end{align}
If we replace $\boldsymbol{\hat{\theta}}_{vr}$ with  $\boldsymbol{\hat{\theta}}{'}_{vr}$, then all the conclusions above hold, including the case that the samples on $\mathcal{I}_1$ are anomalous.
\end{theorem}
\begin{remark}From (\ref{t15}), we have for each entry in $\boldsymbol{\hat{\theta}}_{vr}$ or $\boldsymbol{\hat{\theta}}{'}_{vr}$, the asymptotic  efficiency relative to the mean estimator is $1/D_K$.
	When $K$ goes to infinity,  $D_K$ converges to $\pi/3$ (see Example 1 in \cite{2008Composite} for details), which corresponds to the asymptotic efficiency $0.955$ relative to the mean estimator. For a fixed $K$, we can directly calculate the value of $D_K$ and its asymptotic relative efficiency. For different $K$, the values of $D_K$ and the asymptotic efficiencies relative to the mean estimator are shown in Table \ref{T1}. We can observe that when $K=5$, the asymptotic relative efficiency has exceeded 0.9. When K is greater than 10, the improvement of asymptotic relative efficiency is trivial as $K$ increases. Hence, to balance the computational complexity and the statistical efficiency of the estimator, we can choose $K$ between $5$ and $10$.
	\begin{table}
		\renewcommand\arraystretch{1.0}
			\centering
		\caption{The asymptotic efficiency for different $K$. }
		\label{T1}
		\scalebox{1.0}{
			\begin{tabular}{lccccccccc}
				\hline \quad\quad\quad $K$  	& K=3 & K=5 & K=7 &  K=10 & K=15 & K=20 &K=30 &K=50 &K=100\\
				\hline The value of $D_K$ & 1.168 & 1.103 & 1.080
				& 1.066 &1.056 & 1.053 &1.050 &1.048 &1.047\\
				Asymptotic 
				efficiency & 0.856 & 0.906 & 0.926
				& 0.938 & 0.947 & 0.950 &0.952 &0.954 &0.955\\
				
			    \hline
			\end{tabular}
		}
\end{table}
\end{remark}
\begin{remark}
Together with the results in Lemma \ref{oMom} and Theorem \ref{key1}, we can clearly see that $\|\boldsymbol{\hat{\theta}}_{vr}-\boldsymbol{\theta^*}\|$ has a convergence rate with an improvement of ``$\sqrt{\log n}$'' compared with $\|\boldsymbol{\hat{\theta}}_{med}-\boldsymbol{\theta^*}\|$. When the rate constraints are satisfied, each entry of $\boldsymbol{\hat{\theta}}_{vr}$ has asymptotic normality with a variance that is $D_K$ times that of the global estimator. However, although the asymptotic normality of VRMOL estimator $\boldsymbol{\hat{\theta}}_{vr}$ is provided in (\ref{t14}), it needs very strict rate constraints ($p^2\log^3 n/m=o(1)$ and ${p^3m\log^2 n}/n=o(1)$), which makes its application highly limited when $p$ is divergent. Therefore, we will propose the Robust One-Step Estimator (ROSE) to relax the rate constraints and provide its asymptotic normality in Sections \ref{sec3} and \ref{sec4}.
\end{remark}
\begin{remark}
   Although the expression for $\hat{\theta}_{vr,l}^K$ is much more complex than $\hat{\theta}_{med,l}$, $\hat{\theta}_{vr,l}^K$ has the same order of computational complexity as $\hat{\theta}_{med,l}$. Concretely speaking, for both $\hat{\theta}_{med,l}$ and $\hat{\theta}_{vr,l}^K$, both methods take $O(n)$ time complexity for the node machines to gain the local estimators.  After receiving the local estimators, it takes the central processor $O(m)$ operations to gain $\hat{\theta}_{med,l}$. In contrast, for $\hat{\theta}_{vr,l}^K$, it requires a total of $O(m+n+K)$ time complexity (see Equation (7) in \cite{tu2021variance} for details). Thus, when $K$ is fixed, both $\hat{\theta}_{med,l}$ and $\hat{\theta}_{vr,l}^K$ have the time complexity $O(m+n)$.
\end{remark}

\section{The ROSE algorithm}\label{sec3}
In this section, 
we first restate  the classical one-step estimator (\cite{huo2019distributed}), and then propose a new Robust One-Step Estimator (ROSE) based on the VRMOL estimators.

When all machines are normal, we can calculate the average values of the local gradient  and the local Hessian matrix estimators on node machines with the parameter value $\boldsymbol{\theta}=\boldsymbol{\hat{\theta}}_{vr}$ or $\boldsymbol{\theta}=\boldsymbol{\hat{\theta}}{'}_{vr}$ as the estimators of the gradient and the Hessian matrix.
  According to \cite{huo2019distributed}, the target parameter $\boldsymbol{\theta}$ can be estimated as:
\begin{align}\label{os1}
	\boldsymbol{\hat{\theta}}_{os}=\frac{1}{m}\sum_{j=1}^{m}\boldsymbol{\hat{\theta}}_j-\left[\frac{1}{m}\sum_{j=1}^{m}\nabla^2F_j\left(\frac{1}{m}\sum_{j=1}^{m}\boldsymbol{\hat{\theta}}_j\right)\right]^{-1}\frac{1}{m}\sum_{j=1}^{m}\nabla F_j\left(\frac{1}{m}\sum_{j=1}^{m}\boldsymbol{\hat{\theta}}_j\right).
\end{align}
Under certain regularity assumptions, it holds that for any constant vector $\boldsymbol{v}\in\mathbb{R}^p$ with $\|\boldsymbol{v}\|=1$,
\begin{align*}
	\frac{\sqrt{mn}}{\boldsymbol{v}^{\top}\boldsymbol{\Sigma}(\boldsymbol{\theta^*})\boldsymbol{v}}
	\langle \boldsymbol{\hat{\theta}}_{os}-\boldsymbol{\theta^*},\boldsymbol{v}\rangle\stackrel{d}{\to}\mathbf{N}(0,1).
\end{align*}

\subsection{The VRMOL estimation for the gradient and Hessian matrix}

When Byzantine machines exist, we can use the median of the local estimators to estimate the gradient and the Hessian matrix. The forms of these local estimators are concise: for every $l\in[p]$, define
\begin{align}\label{ag0}
	 \nabla F_{j,l}(\boldsymbol{\theta})=\left\{\begin{array}{cl}
		\nabla_{\theta_l} F_j(\boldsymbol{\theta}) & j \notin \mathcal{B}, \\
		* & j \in \mathcal{B},
	\end{array}
	\right.
	\ \text{and}\
	\nabla^2 F_{j,l_1l_2}(\boldsymbol{\theta})=\left\{\begin{array}{cl}
		\nabla^2_{\theta_{l_1}\theta_{l_2}} F_j(\boldsymbol{\theta}) & j \notin \mathcal{B}, \\
		* & j \in \mathcal{B},
	\end{array}
	\right.
\end{align}
 where $*$ can be any value given by the Byzantine machine.

 	Let
 $\widehat{\nabla F}_{med,l}(\boldsymbol{\theta})=\operatorname{med}\{\nabla_{\theta_l} F_j(\boldsymbol{\theta}),j\in[m]\}$, $l\in[p]$, $\widehat{\nabla F}_{med}(\boldsymbol{\theta})=(\widehat{\nabla F}_{med,1}(\boldsymbol{\theta}),\cdots,\widehat{\nabla F}_{med,p}(\boldsymbol{\theta}))$, $\widehat{\nabla^2 F}_{med,l_1l_2}(\boldsymbol{\theta})=\operatorname{med}\{\nabla^2_{\theta_{l_1}\theta_{l_2}} F_j(\boldsymbol{\theta}),j\in[m]\}$, $l_1,l_2\in[p]$ and $\widehat{\nabla^2 F}_{med}(\boldsymbol{\theta})$ be a matrix with the $(l_1,l_2)$-th entry $\widehat{\nabla^2 F}_{med,l_1l_2}(\boldsymbol{\theta})$.  Then
$\widehat{\nabla F}_{med}(\boldsymbol{\hat{\theta}}_{med})$ and $\widehat{\nabla^2 F}_{med}(\boldsymbol{\hat{\theta}}_{med})$ can be used as the robust estimators of the gradient and the Hessian matrix.
 But as commented before, the median has low asymptotic efficiency. We now use VRMOL to estimate each entry in the gradient and the Hessian matrix.
 It is noteworthy that both the gradient and  Hessian matrix have  simple sample average forms.  The following two methods result in the same estimators of the gradient and the Hessian matrix: the method that we directly use all samples to calculate the gradient and the Hessian matrix at $\boldsymbol{\theta}$, and the method that we first calculate the gradient and the Hessian matrix at $\boldsymbol{\theta}$ on each machine, and then average the results calculated by these machines. In this case, the VRMOL estimators for the gradient and the Hessian matrix have simpler forms than the general M-estimation.

 To be precise,  letting
${\sigma}_{gr,l}^2(\boldsymbol{\theta})=\mathbb{E}\left[\left\{\nabla_{\theta_l}f(\boldsymbol{X}_i,\boldsymbol{\theta})-\nabla_{\theta_l}F_{\mu}(\boldsymbol{\theta})\right\}^2\right] $ and
${\sigma}_{he,l_1l_2}^2(\boldsymbol{\theta})=\mathbb{E}[\{\nabla^2_{\theta_{l_1}\theta_{l_2}}f(\boldsymbol{X}_i,\boldsymbol{\theta})-\nabla^2_{\theta_{l_1}\theta_{l_2}}F_1(\boldsymbol{\theta})\}^2],$ their estimators are defined as
\begin{align}\label{ag2} \hat{\sigma}_{gr,jl}^2(\boldsymbol{\theta})=\left\{\begin{array}{cl}\frac{1}{n}\sum_{i\in\mathcal{I}_j}\{\nabla_{\theta_l}f(\boldsymbol{X}_i,\boldsymbol{\theta})-\nabla_{\theta_l}F_j(\boldsymbol{\theta})\}^2
		& j \notin \mathcal{B}, \\
		* & j \in \mathcal{B},
	\end{array}
	\right.
\end{align}
 and
\begin{align}\label{ag3} \hat{\sigma}_{he,jl_1l_2}^2(\boldsymbol{\theta})=\left\{\begin{array}{cl}\frac{1}{n}\sum_{i\in\mathcal{I}_j}\{\nabla^2_{\theta_{l_1}\theta_{l_2}}f(\boldsymbol{X}_i,\boldsymbol{\theta})-\nabla^2_{\theta_{l_1}\theta_{l_2}}F_j(\boldsymbol{\theta})\}^2
		& j \notin \mathcal{B}, \\
		* & j \in \mathcal{B},
	\end{array}
	\right.
\end{align}
 where $*$ can be  arbitrary values.

Similar to the estimation of $\sigma_l(\boldsymbol{\theta^*})$, if we have the prior information that $\mathcal{B}\subseteq\{2,3,\cdots,m\}$,  we can select $\hat{\sigma}_{gr,1l}^2(\boldsymbol{\theta})$ and $\hat{\sigma}_{he,1l_1l_2}^2(\boldsymbol{\theta})$ as the estimators of ${\sigma}_{gr,l}^2(\boldsymbol{\theta})$ and ${\sigma}_{he,l_1l_2}^2(\boldsymbol{\theta})$ respectively. Otherwise, we use $\hat{\sigma}_{gr,l}^2(\boldsymbol{\theta})=\operatorname{med}\{\hat{\sigma}_{gr,jl}^2(\boldsymbol{\theta}),j\in[m]\}$ and $\hat{\sigma}_{he,l_1l_2}^2(\boldsymbol{\theta})=\operatorname{med}\{\hat{\sigma}_{he,jl_1l_2}^2(\boldsymbol{\theta}),j\in[m]\}$. Further, in the case where $\mathcal{B}\subseteq\{2,3,\cdots,m\}$,  define the $l$-th entry  of $\widehat{\nabla F}_{vr}(\boldsymbol{\theta})=(\widehat{\nabla F}_{vr,1}(\boldsymbol{\theta}),\cdots,\widehat{\nabla F}_{vr,p}(\boldsymbol{\theta}))$ as
\begin{align}\label{0318}
	\widehat{\nabla F}_{vr,l}(\boldsymbol{\theta})=&\widehat{\nabla F}_{med,l}(\boldsymbol{\theta})\nonumber\\
&-\frac{\hat{\sigma}_{gr,1l}(\boldsymbol{\theta})\sum_{k=1}^{K}\sum_{j=1}^{m}\left[\mathbb{I}(\nabla_{\theta_l} F_j(\boldsymbol{\theta})\leq \widehat{\nabla F}_{med,l}(\boldsymbol{\theta})+\hat{\sigma}_{gr,1l}(\boldsymbol{\theta})\Delta_k/\sqrt{n})-\frac{k}{K+1}\right]}{m\sqrt{n}\sum_{k=1}^{K}\psi(\Delta_k)},
\end{align}
and the $(l_1,l_2)$-th entry of the $p\times p$ matrix $\widehat{\nabla^2 F}_{vr}(\boldsymbol{\theta})$ as
 \begin{align}\label{0319}
	&\widehat{\nabla^2 F}_{vr,l_1l_2}(\boldsymbol{\theta})=\widehat{\nabla^2 F}_{med,l_1l_2}(\boldsymbol{\theta})\nonumber
	\\&\ \ \ -\frac{\hat{\sigma}_{he,1l_1l_2}(\boldsymbol{\theta})\sum_{k=1}^{K}\sum_{j=1}^{m}\left[\mathbb{I}(\nabla^2_{\theta_{l_1}\theta_{l_2}} F_j(\boldsymbol{\theta})\leq \widehat{\nabla^2 F}_{med,l_1l_2}(\boldsymbol{\theta})+\hat{\sigma}_{he,1l_1l_2}(\boldsymbol{\theta})\Delta_k/\sqrt{n})-\frac{k}{K+1}\right]}{m\sqrt{n}\sum_{k=1}^{K}\psi(\Delta_k)}.
\end{align}
Then $\widehat{\nabla F}_{vr}(\boldsymbol{\theta})$ and $\widehat{\nabla^2 F}_{vr}(\boldsymbol{\theta})$ are the VRMOL estimators for the gradient and the Hessian matrix at the parameter value $\boldsymbol{\theta}$ seperately.

Therefore,  we can define $\widehat{\nabla F}{'}_{vr,l}(\boldsymbol{\theta})$ and $\widehat{\nabla^2 F}{'}_{vr,l_1l_2}(\boldsymbol{\theta})$ just by replacing $\hat{\sigma}_{gr,1l}(\boldsymbol{\theta})$ and $\hat{\sigma}_{he,1l_1l_2}(\boldsymbol{\theta})$ with $\hat{\sigma}_{gr,l}(\boldsymbol{\theta})$ and $\hat{\sigma}_{he,l_1l_2}(\boldsymbol{\theta})$ respectively, which are still robust when the central processor has anomalous samples.

For the Newton-Raphson iteration, the accuracy of the gradient estimator can be assessed by the difference between the estimated value of the gradient at
$\boldsymbol{\theta}$ and its corresponding expectation at
$\boldsymbol{\theta}$. Analogously, we also concern  the precision of the Hessian matrix estimator. As such, we introduce the following two lemmas.

 \begin{lemma}\label{gh1}
 	Let $\mathcal{B}\subseteq\{2,3,\cdots,n\}$. Suppose that Assumptions \ref{a5} and \ref{a10} hold, $\log^3 n/m=o(1)$, $m\log n/n=o(1)$, $p^2/n=o(1)$ and $\alpha_n=O(1/\log n)$. Then there exists a positive constant $\delta$ such that for $\boldsymbol{\theta}\in B(\boldsymbol{\theta^*},\delta)$,
 	\begin{align*}
 		\|\widehat{\nabla F}_{vr}(\boldsymbol{\theta})-\nabla F_{\mu}(\boldsymbol{\theta})\|=O_p\left(\frac{\alpha_n\sqrt{p}}{\sqrt{n}}
 		+\sqrt{\frac{p}{mn}}\right).
 	\end{align*}
\end{lemma}
\begin{lemma}\label{gh2}
	Let $\mathcal{B}\subseteq\{2,3,\cdots,m\}$. Suppose that Assumptions \ref{a7} and \ref{a102} hold, $\log^3 n/m=o(1)$, $p^2/n=o(1)$, $m\log n/n=o(1)$ and $\alpha_n=o(1/\log n)$.  Then there exists a positive constant $\delta$ such that for $\boldsymbol{\theta}\in B(\boldsymbol{\theta^*},\delta)$,{
	\begin{align*}
		\|\widehat{\nabla^2 F}_{vr}(\boldsymbol{\theta})-\nabla^2 F_{\mu}(\boldsymbol{\theta})\|=O_p\left(\frac{\alpha_np}{\sqrt{n}}
		+\frac{p}{\sqrt{mn}}\right).
	\end{align*}}
\end{lemma}
From these two lemmas, we observe that when $\alpha_n=O(1/\sqrt{m})$, the convergence rates of these two estimators can reach $\sqrt{p/mn}$ and $p/\sqrt{mn}$ respectively, which are the same as the rates if we put all the data on a normal machine.  If we replace $\widehat{\nabla F}_{vr}(\boldsymbol{\theta})$ and $\widehat{\nabla^2 F}_{vr}(\boldsymbol{\theta})$ by $\widehat{\nabla F}{'}_{vr}(\boldsymbol{\theta})$ and $\widehat{\nabla^2 F}{'}_{vr}(\boldsymbol{\theta})$, then when $\mathcal{B}\subseteq\{1,2,3,\cdots,m\}$, the results in these two lemmas also hold. See Remark 1 in Supplementary Material for details.

\subsection{The robust one-step estimator}
Let $r_n=o(1)$ be some rate. For any $\boldsymbol{\theta}$ in $B(\boldsymbol{\theta^*},r_n)$, similar to ${\boldsymbol{\hat{\theta}}}_{os}$ in (\ref{os1}), define two robust one-step estimators as
\begin{align}\label{rose}
	\boldsymbol{\hat{\theta}}_{ros}=\boldsymbol{\theta}-\{\widehat{\nabla^2 F}_{vr}(\boldsymbol{\theta})\}^{-1}\widehat{\nabla F}_{vr}(\boldsymbol{\theta}) \quad\text{and}\quad \boldsymbol{\hat{\theta}}{'}_{ros}=\boldsymbol{\theta}-\{\widehat{\nabla^2 F}{'}_{vr}(\boldsymbol{\theta})\}^{-1}\widehat{\nabla F}{'}_{vr}(\boldsymbol{\theta}).
\end{align}
These estimators are suitable in the case with $\mathcal{B}\subseteq\{2,3,\cdots,m\}$, but $\boldsymbol{\hat{\theta}}{'}_{ros}$ still works when the central processor is a Byzantine machine.
In our research, we mainly analyze the case when the initial value for iteration lies in $B(\boldsymbol{\theta^*},r_n)$ in probability, where $r_n=O({\alpha_n\sqrt{p/n}}+\sqrt{{p}/{mn}}+{p^{3/2}\log n}/{n})$, since both $\|\boldsymbol{\hat{\theta}}_{vr}-\boldsymbol{\theta^*}\|$ and $\|\boldsymbol{\hat{\theta}}{'}_{vr}-\boldsymbol{\theta^*}\|$ are $O_p({\alpha_n\sqrt{p/n}}+\sqrt{{p}/{mn}}+{p^{3/2}\log n}/{n})$.
From now on, unless otherwise specified, we will use $\boldsymbol{\hat{\theta}}_{vr}$($\boldsymbol{\hat{\theta}}{'}_{vr}$) as the initial estimator in the construction of $\boldsymbol{\hat{\theta}}_{ros}$($\boldsymbol{\hat{\theta}}{'}_{ros}$), that is,
\begin{align}
\boldsymbol{\hat{\theta}}_{ros}&=\boldsymbol{\hat{\theta}}_{vr}-\{\widehat{\nabla^2 F}_{vr}(\boldsymbol{\hat{\theta}}_{vr})\}^{-1}\widehat{\nabla F}_{vr}(\boldsymbol{\hat{\theta}}_{vr}),\label{lastros1}\\
\boldsymbol{\hat{\theta}}{'}_{ros}&=\boldsymbol{\hat{\theta}}{'}_{vr}-\{\widehat{\nabla^2 F}{'}_{vr}(\boldsymbol{\hat{\theta}}{'}_{vr})\}^{-1}\widehat{\nabla F}{'}_{vr}(\boldsymbol{\hat{\theta}}{'}_{vr}). \label{lastros2}
\end{align}

Algorithm \ref{alg1} (under normal conditions) and Algorithm \ref{alg2} (allowing for anomalous samples on
$\mathcal{I}_1$) depict the computation procedures involved. Two rounds of communication between the central processor and node machines are required by the proposed algorithms. For Algorithm \ref{alg1}, in the first round of communication, the node machines transmit local parameter estimators to the central processor. The central processor then computes the VRMOL parameter estimator and sends it back to the node machines. In the second round of communication, node machines utilize the VRMOL parameter estimator provided by the central processor to compute  gradients and Hessian matrices, transmitting them back to the central processor. Subsequently, the central processor computes the VRMOL estimators of gradient and Hessian matrix, and computes a one-step correction estimator based on the initial parameter value utilizing the Newton-Raphson iteration. For Algorithm \ref{alg2}, 
the node machines need to transfer several additional variance estimators to the central processor (see the table in Algorithm \ref{alg2} for details). 
This will increase the computational load and information transmission cost of the node machines, but it can handle the cases where  the samples on the central processor are  anomalous or missing.


\begin{remark}We now discuss the order of the time complexity of $\boldsymbol{\hat{\theta}}_{ros}$ and $\boldsymbol{\hat{\theta}}{'}_{ros}$.
	When $K$ is fixed,  $\boldsymbol{\hat{\theta}}_{vr}$, $\boldsymbol{\hat{\theta}}{'}_{vr}$, $\widehat{\nabla F}_{vr}(\boldsymbol{\hat{\theta}}_{vr})$ and $\widehat{\nabla F}{'}_{vr}(\boldsymbol{\hat{\theta}}{'}_{vr})$  take $O(p(m+n))$ time complexity, and $\widehat{\nabla^2 F}_{vr}(\boldsymbol{\hat{\theta}}_{vr})$ and $\widehat{\nabla^2 F}{'}_{vr}(\boldsymbol{\hat{\theta}}{'}_{vr})$  take $O(p^2(m+n))$ time complexity. So both $\boldsymbol{\hat{\theta}}_{ros}$ and $\boldsymbol{\hat{\theta}}{'}_{ros}$ have the total time complexity $O(p^2(m+n))$. As these machines can work simultaneously, there is no significant difference in computational time between the  variance estimators $\{\hat{\sigma}_l(\boldsymbol{\hat{\theta}}_{med}),\hat{\sigma}_{gr,1l}(\boldsymbol{\hat{\theta}}_{vr}),\hat{\sigma}_{he,1l_1l_2}(\boldsymbol{\hat{\theta}}_{vr})\}$ and  $\{\hat{\sigma}{'}_l,\hat{\sigma}_{gr,l}(\boldsymbol{\hat{\theta}}{'}_{vr}),\hat{\sigma}_{he,l_1l_2}(\boldsymbol{\hat{\theta}}{'}_{vr})\}$. However, as shown in the algorithms, the amount of data transmitted from the node machines  to the central processor of Algorithm \ref{alg2} is twice  of that in Algorithm \ref{alg1}, so the communication  cost of Algorithm \ref{alg2} is higher. Therefore, if we have prior information that the central processor is normal, we can prioritize using Algorithm \ref{alg1}.
\end{remark}
\begin{algorithm}
	\DontPrintSemicolon
	\KwInput{{The dataset $\{X_1,X_2,\dots,X_N\}$ which is evenly distributed on $m$ machines $\{\mathcal{I}_1,\mathcal{I}_2,\dots,\mathcal{I}_m\}$ with the local sample size $n$. A  positive integer $K$.}}

	Each machine computes a local M-estimator ${\boldsymbol{\hat{\theta}}}_j$, $j\in[m]$ by minimizing the local loss function $F_j(\boldsymbol{\theta})$ and then sends ${\boldsymbol{\hat{\theta}}}_j$ to the central processor $\mathcal{I}_1$.
	
	Based on the central processor $\mathcal{I}_1$,  calculate $\boldsymbol{\hat{\theta}}_{med}=(\hat{\theta}_{med,1},\cdots,\hat{\theta}_{med,p})^{\t}$ and $\hat{\sigma}_l(\boldsymbol{\hat{\theta}}_{med})$  in the formula (\ref{sig2}), where
	$\hat{\theta}_{med,l}=\operatorname{med}\{\hat{\theta}_{jl},j\in[m]\}$, $l\in[p]$, $\boldsymbol{\hat{\theta}}_{vr}=(\hat{\theta}_{vr,1}^K,\hat{\theta}_{vr,2}^K,\cdots,\hat{\theta}_{vr,p}^K)^{\top}$ with $\hat{\theta}_{vr,l}^K$ estimated based on VRMOL in  (\ref{vr4}).

	Distribute $\boldsymbol{\hat{\theta}}_{vr}$ to each machine $\mathcal{I}_j, j\in[m]$.
	
	{For $j\in[m]$, the $j$-th machine computes the local gradient $\nabla F_{j,l}(\boldsymbol{\hat{\theta}}_{vr})$, $l\in[p]$ and the local Hessian matrix  $\nabla^2 F_{j,l_1l_2}(\boldsymbol{\hat{\theta}}_{vr})$, $l_1,l_2\in[p]$ as (\ref{ag0})}.

	Then the $j$-th machine sends
	$\nabla F_{j,l}(\boldsymbol{\hat{\theta}}_{vr})$ and
	$\nabla^2 F_{j,l_1l_2}(\boldsymbol{\hat{\theta}}_{vr})$ back to the central processor $\mathcal{I}_1$.

	The central processor $\mathcal{I}_1$ constructs the VRMOL aggregated gradient $\widehat{\nabla F}_{vr}(\boldsymbol{\hat{\theta}}_{vr})=(\widehat{\nabla F}_{vr,1}(\boldsymbol{\hat{\theta}}_{vr}),\cdots,\widehat{\nabla F}_{vr,p}(\boldsymbol{\hat{\theta}}_{vr}))$,
	where each $l$-th coordinate is calculated as (\ref{0318}),
 and the VRMOL aggregated Hessian matrix $	\widehat{\nabla^2 F}_{vr}(\boldsymbol{\hat{\theta}}_{vr})$ with the $(l_1,l_2)$-th entry calculated by (\ref{0319}).

	 Calculates $\boldsymbol{\hat{\theta}}_{ros}$ based on the central processor $\mathcal{I}_1$ as the formula (\ref{lastros1}).

	\KwOutput{The final estimator $\boldsymbol{\hat{\theta}}_{ros}$.}
	\caption{The ROSE $\boldsymbol{\hat{\theta}}_{ros}$ (Normal central processor)}
	\label{alg1}
\end{algorithm}

\begin{algorithm}
	\DontPrintSemicolon
	\KwInput{The same as Algorithm \ref{alg1}.}
Each machine computes  ${\boldsymbol{\hat{\theta}}}_j=(\hat{\theta}_{j1},\cdots,\hat{\theta}_{jp})$ and $\hat{\sigma}_{jl}(\boldsymbol{\hat{\theta}}_j)$,  $j\in[m]$ and then sends them to the central processor $\mathcal{I}_1$.
	
	For $l\in[p]$, based on the central processor $\mathcal{I}_1$, calculate $\hat{\theta}_{med,l}=\operatorname{med}\{\hat{\theta}_{jl},j\in[m]\}$, $\hat{\sigma}{'}_{l}=\operatorname{med}\{\hat{\sigma}_{jl}(\boldsymbol{\hat{\theta}}_j), j\in[m]\}$, and $\hat{\theta}{'}_{vr,l}^K$ by VRMOL in  (\ref{vr4})(replacing $\hat{\sigma}_l(\boldsymbol{\hat{\theta}}_{med})$ by $\hat{\sigma}{'}_{l}$).

	Distribute $\boldsymbol{\hat{\theta}}{'}_{vr}=(\hat{\theta}{'}_{vr,1}^K,\hat{\theta}{'}_{vr,2}^K,\cdots,\hat{\theta}{'}_{vr,p}^K)$ to each machine $\mathcal{I}_j, j\in[m]$.

	For the $j$-th machine, compute  $\nabla F_{j,l}(\boldsymbol{\hat{\theta}}{'}_{vr})$, $l\in[p]$ and  $\nabla^2 F_{j,l_1l_2}(\boldsymbol{\hat{\theta}}{'}_{vr})$, $l_1,l_2\in[p]$ according to (\ref{ag0}), and the corresponding variance estimator $\hat{\sigma}_{gr,jl}^2(\boldsymbol{\hat{\theta}}{'}_{vr})$ and $\hat{\sigma}_{he,jl_1l_2}^2(\boldsymbol{\hat{\theta}}{'}_{vr})$ based on (\ref{ag2}) and (\ref{ag3}) respectively, and then send these results to the central processor $\mathcal{I}_1$.

The central processor $\mathcal{I}_1$ constructs the VRMOL aggregated gradient $\widehat{\nabla F}{'}_{vr}(\boldsymbol{\hat{\theta}}{'}_{vr})$ with each coordinate calculated by the formula (\ref{0318}), but replacing $\hat{\sigma}_{gr,1l}(\boldsymbol{\hat{\theta}}_{vr})$ by $\hat{\sigma}_{gr,l}^2(\boldsymbol{\hat{\theta}}{'}_{vr})$,
	and the VRMOL aggregated Hessian matrix $	\widehat{\nabla^2 F}{'}_{vr}(\boldsymbol{\hat{\theta}}{'}_{vr})$ with the $(l_1,l_2)$-th entry calculated based on the formula (\ref{0319}), but replacing $\hat{\sigma}_{he,1l_1l_2}(\boldsymbol{\hat{\theta}}_{vr})$ by  $\hat{\sigma}_{he,l_1l_2}^2(\boldsymbol{\hat{\theta}}{'}_{vr})$.

 Calculates $\boldsymbol{\hat{\theta}}{'}_{ros}$ based on the central processor $\mathcal{I}_1$ as the formula (\ref{lastros2}).

\KwOutput{The final estimator $\boldsymbol{\hat{\theta}}{'}_{ros}$.}
	\caption{ROSE $\boldsymbol{\hat{\theta}}{'}_{ros}$ (Byzantine central processor)}
	\label{alg2}
\end{algorithm}
\section{Theoretical Results for  ROSE}\label{sec4}
In this section, we present the asymptotic normality and convergence rate of the estimators $\boldsymbol{\hat{\theta}}_{ros}$ and $\boldsymbol{\hat{\theta}}{'}_{ros}$.

\subsection{Asymptotic normality}
 The asymptotic normality of the estimators is stated below.
   \begin{theorem}\label{the2}Suppose Assumptions \ref{a3}, \ref{a5}-\ref{a7} and \ref{a10}-\ref{a12} hold. Let $\alpha_n=o(1/\sqrt{pm}), pm\log n/n=o(1)$, $p\log^2 n/m=o(1)$,  $p^2\sqrt{m}\log^{3/2} n/n=o(1)$ and $r_n=O(\alpha_n\sqrt{p/n}+\sqrt{p/(mn)}+p^{3/2}\log n/n)$. For any constant vector $\boldsymbol{v}$ satisfying $\|\boldsymbol{v}\|=1$, $\boldsymbol{\theta}\in B(\boldsymbol{\theta^*},r_n)$,
  	we have
  	\begin{align*}
  		\frac{\sqrt{mn}}{\sqrt{\boldsymbol{v}^{\top}\mathbf{\Sigma}_{vr}^{K}(\boldsymbol{\theta^*})\boldsymbol{v}}}
  		\langle \boldsymbol{\hat{\theta}}_{ros}-\boldsymbol{\theta^*},\boldsymbol{v}\rangle\stackrel{d}{\to}\mathbf{N}(0,1)\quad\text{and}\quad \frac{\sqrt{mn}}{\sqrt{\boldsymbol{v}^{\top}\mathbf{\Sigma}_{vr}^K(\boldsymbol{\theta^*})\boldsymbol{v}}}
  		\langle \boldsymbol{\hat{\theta}}{'}_{ros}-\boldsymbol{\theta^*},\boldsymbol{v}\rangle\stackrel{d}{\to}\mathbf{N}(0,1),
  	\end{align*}
  where $\boldsymbol{\hat{\theta}}_{ros}$ and $\boldsymbol{\hat{\theta}}{'}_{ros}$ are defined in (\ref{rose}), and the asymptotic normality of $\boldsymbol{\hat{\theta}}_{ros}$ only holds when $\mathcal{B}\subseteq \{2,3,\cdots,m\}$.
  \end{theorem}
 \begin{remark}
 	In Theorem \ref{the2}, we require the initial estimator before iteration has the same convergence rate of $\boldsymbol{\hat{\theta}}_{vr}$ and $\boldsymbol{\hat{\theta}}{'}_{vr}$ proposed in Section \ref{sec2}.  That is to say, if we use another robust estimator $\boldsymbol{\hat{\theta}}_{robust}$ satisfying $\|\boldsymbol{\hat{\theta}}_{robust}-\boldsymbol{\theta^*}\|=O_p({\alpha_n\sqrt{p/n}}+\sqrt{{p}/{mn}}+{p^{3/2}\log n}/{n})$ as the initial value rather than $\boldsymbol{\hat{\theta}}_{vr}$ or $\boldsymbol{\hat{\theta}}{'}_{vr}$, we can gain the same limit normal distribution under the same assumptions. This implies Theorem \ref{the2} has a wider range of applications.
 \end{remark}

Although the  limitations on the relationship between $m$, $n$, and $p$ in Theorem \ref{the2} exist, they  are reasonably much weaker than those for  the asymptotic normality of the initial value $\boldsymbol{\hat{\theta}}_{vr}$ or $\boldsymbol{\hat{\theta}}{'}_{vr}$, as shown in Theorem \ref{key1}. We discuss how the rate constraints for $m$, $n$, and $p$ in Theorem \ref{the2} can be satisfied. For example, if $m\asymp p\log^{3} n$ and $n\asymp p^{5/2}\log^5 n$, then all three rate constraints $pm\log n/n=o(1)$, $p\log^2 n/m=o(1)$ and $p^2\sqrt{m}\log^{3/2} n/n=o(1)$ hold.  Thus the leading terms of $m$ and $n$ are of the orders  $ p$ and $ p^{5/2}$ respectively. It is worth noting that these rate requirements are weaker than Theorem 20 in \cite{tu2021variance}, which requires $p^2\log^3 n/m=o(1)$ and $pm/n=o(1)$ to gain the asymptotic normality. As for the rate connection between $p$ and the total sample size $N=mn$, ignoring the logarithmic term,  \cite{tu2021variance} need $p^5\asymp N$ at least, while our estimator $\boldsymbol{\hat{\theta}}_{ros}$ or $\boldsymbol{\hat{\theta}}{'}_{ros}$ only needs $p^{7/2}\asymp N$. As a comparison, the asymptotic normality of $\boldsymbol{\hat{\theta}}_{vr}$ or $\boldsymbol{\hat{\theta}}{'}_{vr}$ requires at least $p^7\asymp N$, therefore, when $p$ diverges, a one-step correction of $\boldsymbol{\hat{\theta}}_{vr}(\boldsymbol{\hat{\theta}}{'}_{vr})$ is necessary.

As an application of Theorem \ref{the2}, for any constant vector $\boldsymbol{v}$, we can get two $100(1-\alpha)$\% confidence intervals of $\boldsymbol{v}^{\top}\boldsymbol{\theta^*}$ as:
\begin{align*} \mathbf{CI1}_{1-\alpha}(\boldsymbol{v}^{\top}\boldsymbol{\theta^*})&=\left[\boldsymbol{v}^{\top}\boldsymbol{\hat{\theta}}_{ros}-z_{1-\alpha/2}\sqrt{\frac{\boldsymbol{v}^{\top}\widehat{\boldsymbol{\Sigma}_{vr}^K}(\boldsymbol{\theta^*})\boldsymbol{v}}{mn}}, \boldsymbol{v}^{\top}\boldsymbol{\hat{\theta}}_{ros}+z_{1-\alpha/2}\sqrt{\frac{\boldsymbol{v}^{\top}\widehat{\boldsymbol{\Sigma}_{vr}^K}(\boldsymbol{\theta^*})\boldsymbol{v}}{mn}}\right],\\
\mathbf{CI2}_{1-\alpha}(\boldsymbol{v}^{\top}\boldsymbol{\theta^*})&=\left[\boldsymbol{v}^{\top}\boldsymbol{\hat{\theta}}'_{ros}-z_{1-\alpha/2}\sqrt{\frac{\boldsymbol{v}^{\top}\widehat{\boldsymbol{\Sigma}_{vr}^K}(\boldsymbol{\theta^*})\boldsymbol{v}}{mn}}, \boldsymbol{v}^{\top}\boldsymbol{\hat{\theta}}'_{ros}+z_{1-\alpha/2}\sqrt{\frac{\boldsymbol{v}^{\top}\widehat{\boldsymbol{\Sigma}_{vr}^K}(\boldsymbol{\theta^*})\boldsymbol{v}}{mn}}\right],
\end{align*}
where $z_{\alpha}$ is the $\alpha$-quantile of the standard normal distribution, and $\widehat{\boldsymbol{\Sigma}_{vr}^K}(\boldsymbol{\theta^*})$ is a consistent estimator of $\boldsymbol{\Sigma}_{vr}^K(\boldsymbol{\theta^*})$.
\begin{remark}
The quantity	$\widehat{\boldsymbol{\Sigma}_{vr}^K}(\boldsymbol{\theta^*})$ can be estimated by many methods. For example, since $\boldsymbol{\Sigma}(\boldsymbol{\theta^*})$ is close to $\boldsymbol{\Sigma}_{vr}^{K}(\boldsymbol{\theta^*})$ and the entries on the diagonal of $\boldsymbol{\Sigma}_{vr}^{K}(\boldsymbol{\theta^*})$ are  $D_K$ times that of          $\boldsymbol{\Sigma}(\boldsymbol{\theta^*})$, we can convert estimating $\boldsymbol{\Sigma}_{vr}^{K}(\boldsymbol{\theta^*})$ into estimating $\boldsymbol{\Sigma}(\boldsymbol{\theta^*})=\{\nabla^2F_{\mu}(\boldsymbol{\theta^*})\}^{-1}\mathbb{E}\{\nabla f(\boldsymbol{X},\boldsymbol{\theta^*})^{\otimes 2}\} \{\nabla^2F_{\mu}(\boldsymbol{\theta^*})\}^{-1}$. Since $\nabla^2F_{\mu}(\boldsymbol{\theta^*})$ and $\mathbb{E}\{\nabla f(\boldsymbol{X},\boldsymbol{\theta^*})^{\otimes 2}\}$ can be estimated by $\widehat{\nabla^2F}_{vr}(\boldsymbol{\hat{\theta}}_{vr})$ and $1/n\sum_{i\in\mathcal{I}_1}\nabla f(\boldsymbol{X}_i,\boldsymbol{\hat{\theta}}_{ros})^{\otimes 2}$ respectively, we can use the plug-in  estimation for $\boldsymbol{\Sigma}(\boldsymbol{\theta^*})$. This  method does not require additional communication between machines. There are other methods for estimating $\boldsymbol{\Sigma}_{vr}^{K}(\boldsymbol{\theta^*})$, but  this is not the focus of our research, thus, we do not discuss the convergence rates of different estimation methods in detail.
\end{remark}
\subsection{Convergence rate}
If we do not require $\boldsymbol{\hat{\theta}}_{ros}(\boldsymbol{\hat{\theta}}{'}_{ros})$ to have the asymptotic normality but only care about the convergence rate of $\boldsymbol{\hat{\theta}}_{ros}(\boldsymbol{\hat{\theta}}{'}_{ros})$, then the rate constraints required for asymptotic normality can be further relaxed.
\begin{theorem}\label{mse}
Suppose that Assumptions \ref{a3}, \ref{a5}-\ref{a7} and \ref{a10}-\ref{a12} hold and $\mathcal{L}_1$ be a normal machine. Further suppose $\alpha_n=o(1/\log n)$, $\log^3 n/m=o(1)$, $m\log n/n=o(1)$ and $p^2\log^2 n/n=o(1)$. Write $r_n=O(\alpha_n\sqrt{p/n}+\sqrt{p/(mn)}+p^{3/2}\log n/n)$. For any $\boldsymbol{\theta}\in B(\boldsymbol{\theta^*},r_n)$ and $$\boldsymbol{\hat{\theta}}_{ros}=\boldsymbol{\theta}-\{\widehat{\nabla^2 F}_{vr}(\boldsymbol{\theta})\}^{-1}\widehat{\nabla F}_{vr}(\boldsymbol{\theta}),$$ then
	\begin{align*}		\boldsymbol{\hat{\theta}}_{ros}-\boldsymbol{\theta^*}=-[\nabla^2F_{\mu}(\boldsymbol{\theta^*})]^{-1}\widehat{\nabla F}_{vr}(\boldsymbol{\theta^*})+\boldsymbol{R}_N.
	\end{align*}
	Here
	\begin{align*}
		\mathbb{E}[\|[\nabla^2F_{\mu}(\boldsymbol{\theta^*})]^{-1}\widehat{\nabla F}_{vr}(\boldsymbol{\theta^*})\|^2]=&\frac{1}{mn}\operatorname{Tr}([\nabla^2F_{\mu}(\boldsymbol{\theta^*})]^{-1}\mathbf
		{V}_{g,vr}(\boldsymbol{\theta^*})[\nabla^2F_{\mu}(\boldsymbol{\theta^*})]^{-1})
		\\&+O\left(\frac{\alpha_n^2p}{n}\right)+O\left(\frac{\alpha_np}{\sqrt{m}n}\right)+o\left(\frac{p}{mn}\right),
	\end{align*}
where $\mathbf
{V}_{g,vr}(\boldsymbol{\theta})$ is a $p\times p$ matrix with the $(l_1,l_2)$-th entry $$\frac{{\sum_{k_1=1}^K\sum_{k_2=1}^K(\kappa_{k_1,k_2,K}^{l_1,l_2}-\kappa_{k_1}\kappa_{k_2})}\sigma_{gr,l_1}(\boldsymbol{\theta})\sigma_{gr,l_2}(\boldsymbol{\theta})}{\{\sum_{k=1}^{K}\psi(\Delta_{k})\}^2}.$$

When $\log^3 n/m=o(1)$ and $m\log n/n=o(1)$,  we have the following results for the remainder  $\boldsymbol{R}_N$.
	
	(1) If $\alpha_n=o(1/\log n)$, $p^2\log^2 n/n=o(1)$,  then $$\|\boldsymbol{R}_N\|=o_p\left({\sqrt{\frac{p}{n}}}\right).$$
	
	(2) If $\alpha_n=o(1/\sqrt{m})$, $p^2\log^2 n/n=o(1)$,  then $$\|\boldsymbol{R}_N\|=o_p\left(\sqrt{\frac{p}{mn}}+\frac{p}{n}\right).$$
	
	(3) If $\alpha_n=o(1/\sqrt{pm})$, $p^{5/2}/n=o(1)$,  then $$\|\boldsymbol{R}_N\|=o_p\left(\frac{1}{\sqrt{mn}}+\frac{1}{n^{4/5}}\right).$$
	
	(4) If $\alpha_n=o(1/\sqrt{m})$, $p$ is fixed,  then $$\|\boldsymbol{R}_N\|=o_p\left(\frac{1}{\sqrt{mn}}\right).$$
	If we replace $\boldsymbol{\hat{\theta}}_{ros}$, 
	$\widehat{\nabla F}_{vr}(\boldsymbol{\theta})$ and $\widehat{\nabla^2 F}_{vr}(\boldsymbol{\theta})$ with $\boldsymbol{\hat{\theta}}{'}_{ros}$, 
	$\widehat{\nabla F}{'}_{vr}(\boldsymbol{\theta})$ and $\widehat{\nabla^2 F}{'}_{vr}(\boldsymbol{\theta})$ respectively, and  the samples on $\mathcal{I}_1$ are anomalous, the  limit properties are the same as the above.
\end{theorem}
For the relationship between $m$, $n$ and $p$, if we ignore the logarithmic terms, Theorem \ref{mse} only need  $p^2=o(n)$ and $m=o(n)$. In general, the rate constraint on the leading term for the consistency of the M-estimator with the dimension $p$ and the sample size $n$ is $p^2=o(n)$, which is weaker than $p^3=o(n)$ in Assumption G of \cite{tu2021variance}. Since we need each machine to  provide a consistent parameter estimator in the construction of the initial estimator $\boldsymbol{\hat{\theta}}_{vr}$, the rate of the leading term  $p^2=o(n)$ is hard to be further relaxed. $m=o(n)$ is also required in \cite{tu2021variance}, and it is necessary when we use the VRMOL estimator to estimate the gradient. We add an additional ``$\log n$" and require $m\log n/n=o(1)$ since this can make the results more concise.

Moreover, from Theorem \ref{mse}, we can infer that if $\alpha_n=o(1/\sqrt{m})$, $\log^3 n/m=o(1)$, $p^2\log^2 n/n=o(1)$ and $\sqrt{p}m\log n/n=o(1)$, then $\|\boldsymbol{\hat{\theta}}_{ros}-\boldsymbol{\theta^*}\|=O_p(\sqrt{p/mn})$, which has the same convergence rate as the estimator obtained by aggregating all data onto one normal machine. This rate is the same as that derived by the one-step estimator in \cite{huo2019distributed} (for fixed $p$) and the optimal rate $O_p(\alpha_n/\sqrt{n}+\sqrt{p/mn})$ proved by \cite{2018Byzantine} (because $\alpha_n=o(1/\sqrt{m})$), and  faster than that of the estimator provided by \cite{tu2021variance}, which cannot converge to $\boldsymbol{\theta^*}$ with a rate faster than $O(\sqrt{p\log n/mn})$. 

\begin{remark}
	We can also estimate the gradient and the Hessian matrix by computing the medians of the gradient and the Hessian matrix estimators from node machines. However, the estimation errors of these two medians will be larger compared with the VRMOL estimator, which ultimately result in greater estimation error of the final parameter estimator. We will demonstrate this in numerical simulations.
\end{remark}
\begin{remark}
	 $\boldsymbol{\hat{\theta}}{'}_{ros}$ could be applied  to the situations where the samples on the central processor are anomalous or not  available whereas  \cite{Jordan2019Communication} and \cite{tu2021variance} require the presence of samples on the central processor. If we generate numerical values at random to be employed as the local statistics computed by the central processor when no actual samples are available, the theoretical properties of $\boldsymbol{\hat{\theta}}{'}_{ros}$ retain the same.

\end{remark}

\section{Numerical studies}\label{sec5}
 We conduct several numerical experiments in this section. For the proposed algorithms,  we make a comparison among  five distinct distributed computing methods using the  average root squared errors in our simulations. For synthetic data, we consider  the logistic   and Poisson regression models. For real data example, we employ the logistic regression model to analyze the MNIST dataset,  which is widely used in machine learning research.

\subsection{Synthetic data}\label{5.1}
In this subsection, we  conduct distributed logistic  and distributed Poisson regressions to examine the effectiveness and robustness of the proposed ROSE estimator.  Subsections \ref{5.1.1} and \ref{5.1.2} consider that the central processor $\mathcal{I}_1$ is normal and compare Algorithm \ref{alg1} with other methods. Subsection \ref{5.1.3} conducts the case that the central processor has anomalous data and demonstrate the advantages of Algorithm \ref{alg2} through comparison.
\subsubsection{The results for logistic regression}\label{5.1.1}

Consider $\mathcal{B}\subseteq \{2,3,\cdots,m\}$ and  a general logistic regression model as:
\begin{align}\label{log2}
	Y\sim \operatorname{Bernoulli}(p_0), \text{where}\quad p_0=\frac{\exp(\boldsymbol{X}^{\top}\boldsymbol{\theta^*})}{1+\exp(\boldsymbol{X}^{\top}\boldsymbol{\theta^*})},
\end{align}
where $Y\in\{0,1\}$ is a binary response variable, $\boldsymbol{X}\in\mathbb{R}^q$ follows the multivariate normal distribution $\mathbf{N}(\boldsymbol{0},\boldsymbol{\Sigma}_T)$ with $\boldsymbol{\Sigma}_T$ being a symmetric Toeplitz matrix with the $(i,j)$-entry $(0.5)^{|i-j|}$, $i,j\in[q]$, and $\boldsymbol{\theta^*}\in\mathbb{R}^p$ is the interesting parameter vector and $\boldsymbol{\theta^*}=p^{-1/2}(1,(p-2)/(p-1),(p-3)/(p-1),\cdots,0)^{\top}$.  We set the dimension $p=q=30$, the number of machines $m=11,31,101$ (including one central processor), and the sample size on each machine $n=200,300,500,1000$ respectively. Then $N=mn$ is the total sample size.

  In this experiment, we generate the Byzantine machines by scaling attacks based on $-3$ times normal value. That is to say, the statistics transmitted by the node machines to the central processor are $-3$ times the correct statistics. 
  Consider three ratios of Byzantine machines ($\alpha_n$): (1) $\alpha_n=0.00$, (2) $\alpha_n=0.10$, and (3) $\alpha_n=0.20$. To illustrate the effectiveness and robustness of  ROSE, we compare the following six estimation methods: (1) the standard averaging
estimator (proposed by \cite{2013Communication}); (2) the one-step estimator (proposed by \cite{huo2019distributed}); (3)  ROSE (med), which uses $\boldsymbol{\hat{\theta}}_{med}$, $\widehat{\nabla F}_{med}(\boldsymbol{\hat{\theta}}_{med})$ and $\widehat{\nabla^2 F}_{med}(\boldsymbol{\hat{\theta}}_{med})$ as the initial value, the gradient and the Hessian matrix respectively in our one-step approach; (4) ROSE (K=10), computed by Algorithm \ref{alg1}, which uses $\boldsymbol{\hat{\theta}}_{vr}$, $\widehat{\nabla F}_{vr}(\boldsymbol{\hat{\theta}}_{vr})$ and $\widehat{\nabla^2 F}_{vr}(\boldsymbol{\hat{\theta}}_{vr})$ as the initial value, the gradient and the Hessian matrix respectively, and compute $\boldsymbol{\hat{\theta}}_{ros}$ with $K=10$; (5) the CSL approximation proposed by \cite{Jordan2019Communication}; and (6) the robust approximation based on Algorithm \ref{alg1} with $K=10$  which was proposed by \cite{tu2021variance}). Here (1), (2), and (5) are non-robust methods, while (3), (4), and (6) are robust. Regarding the number of iterations necessary for the algorithmic convergence, (1) does not require iteration, (2), (3), and (4) necessitate only a single round of iteration, and (5) and (6) demand multiple rounds of iteration. {As  \cite{tu2021variance} used,  for the methods (5) and (6), we iterate ten rounds in the actual simulation.}
To compute the average root squared error and its corresponding standard error of the parameter estimator, we replicate each experiment 100 times. The results generated from these simulations are reported in Table \ref{T2}. For each cell in the table, the first value corresponds to the average root squared error obtained over the 100 experiments, while the number enclosed in parentheses denotes the standard deviation of the root squared error.

Table \ref{T2} reports that ROSE (K=10) consistently outperforms the standard averaging estimator proposed by \cite{2013Communication} and ROSE (med). Although the methods posited by \cite{Jordan2019Communication} and \cite{tu2021variance} are ineffective (the estimators fail to converge)
when the local sample size, $n$, assumes values of $200$ or $300$,   ROSE (K=10) continues to provide robust estimation outcomes. When all machines are normal and  the local sample size remains small, ROSE (K=10) may perform slightly worse than \cite{huo2019distributed}. However, when the local sample size increases to $n=500$ or $1000$, this difference becomes negligible. When 10\% or 20\%  Byzantine machines exist and local sample sizes remain $n=200$ or $300$, ROSE (K=10) yields significantly better results than the estimators proposed by \cite{2013Communication}, \cite{huo2019distributed}, and \cite{Jordan2019Communication} (exclusive of cases where  $m=11$ and $\alpha_n=0.10$). Additionally, the performance of ROSE (K=10) is slightly better that that of \cite{tu2021variance} when $n=500$ or $1000$,  although ROSE (K=10) requires merely one round of iteration.

\subsubsection{The results for Poisson regression}\label{5.1.2}
Consider $\mathcal{B}\subseteq \{2,3,\cdots,m\}$ and  the Poisson regression model as:
\begin{align}\label{poi2}
	Y\sim \operatorname{Poission}(\boldsymbol{\lambda}),\ \text{with}\ \boldsymbol{\lambda}={\exp(\boldsymbol{X}^{\top}\boldsymbol{\theta^*})},
\end{align}
where the response variable $Y$ follows the poisson distribution with the parameter $\boldsymbol{\lambda}$, and all parameter values are kept consistent with those of the logistic regression model. The data is generated following (\ref{poi2}) on normal machines, while the response variable $Y$ on Byzantine machines conform to $\operatorname{Poission}(\exp(-\boldsymbol{X}^{\top}{\boldsymbol{\theta^*}}))$. The average and standard deviation of  root squared errors are shown in Table \ref{T3}. In each cell, the first number corresponds to the average of root squared errors obtained across 100 experiments, while the number in parentheses denotes the associated standard deviation.

Table \ref{T3} indicates that  ROSE (K=10) performs better than the methods introduced by \cite{2013Communication} and  ROSE (med) in most scenarios. When Byzantine machines do not exist, the average of  root squared errors of  ROSE (K=10) is almost identical to that of \cite{huo2019distributed}. On the other hand, when the proportion of Byzantine machines is 0.10 or 0.20,  ROSE (K=10) provides considerably better results than the alternative methods in comparison. For $n=200$, $300$ or $500$, the estimators provided by \cite{Jordan2019Communication} and \cite{tu2021variance} perform worse than other approaches. For $n=1000$,  ROSE (K=10)  slightly outperforms \cite{tu2021variance}  when 10\% or 20\% machines are Byzantine machines, but outperforms \cite{Jordan2019Communication}.

Besides, 
it is noteworthy that the approaches proposed by \cite{Jordan2019Communication} and \cite{tu2021variance} necessitate a relatively large local sample size. These two methods employ Hessian matrix information but only estimate by the data from central processor with a sample size of $n$. When $p$ is large and $n$ is relatively small, the estimation error associated with the Hessian matrix may be large, leading to non-convergence of the gradient descent method. In our simulations, when $n=200$ or $300$, the estimators that are based on these two methods exhibit considerable deviation from the true parameter values after 10 iterations, so we only report the results  after 3 iterations.  Additionally, in the case where $n=500$, the results provided by these two methods after 10 rounds of iteration show significant deviations from the expected parameter values for the Poisson regression model, as shown in Table \ref{T3}.

\begin{table}
	\renewcommand\arraystretch{1.0}
	\centering
	\caption{Logistic regression model. }
	\label{T2}
	\scalebox{0.8}{
		\begin{threeparttable}
	\begin{tabular}{lccccc}
		\hline \quad\quad\quad m=11  	& $\alpha_n$ & $n=200$ & $n=300$ &   $n=500$ & $n=1000$ \\
		\hline \cite{2013Communication} & 0.00 & 0.4602(0.0688)   & 0.3266(0.0439)   & 0.2335(0.0374) & 0.1600(0.0213)  \\
		& 0.10 & 0.5938(0.0957)   & 0.4610(0.0772)   & 0.3590(0.0482) & 0.2784(0.0262)  \\
		& 0.20 & 0.7923(0.1175)   & 0.6473(0.0738)   & 0.5535(0.0449) & 0.4829(0.0286)  \\
		\cite{huo2019distributed} 	& 0.00 & 0.3222(0.0491)   & 0.2663(0.0373)   & 0.2084(0.0342) & 0.1505(0.0205)  \\
		& 0.10 & 0.7691(0.1456)   & 0.5895(0.1089)   & 0.4342(0.0795) & 0.2927(0.0496)  \\
		& 0.20 & 22.285(51.983) & 17.575(29.670) & 2.2512(1.5873) & 0.8934(0.1695)  \\
		ROSE(med) & 0.00 & 0.4991(0.0890)    & 0.3887(0.0691)   & 0.3055(0.0502) & 0.2104(0.0310)   \\
		& 0.10 & 0.5645(0.0998)   & 0.4466(0.0747)   & 0.3352(0.0537) & 0.2355(0.0363)  \\
		& 0.20 & 0.6615(0.1100)   & 0.5081(0.0984)   & 0.3893(0.0632) & 0.2661(0.0386)  \\
		 ROSE(K=10) & 0.00 & 0.3497(0.0572)   & 0.2858(0.0375)   & 0.2208(0.0324) & 0.1574(0.0251)  \\
		& 0.10 & 0.4863(0.0826)   & 0.3861(0.0563)   & 0.2958(0.0490) & 0.2040(0.0328)  \\
		& 0.20 & 0.5873(0.0944)   & 0.4511(0.0699)   & 0.3537(0.0563) & 0.2378(0.0307)  \\
		\cite{Jordan2019Communication} & 0.00 & 21.430(91.715)*    & 1.5500(9.6715)   & 0.2157(0.0332) & 0.1479(0.0239)  \\
		& 0.10 & 0.8309(0.1574)   & 0.6309(0.1049)   & 0.4536(0.0770) & 0.3155(0.0547)  \\
		& 0.20 & 2.4351(0.8360)*  & 1.4306(0.3199)   & 1.9456(0.6206) & 1.0417(0.2243)  \\
		\cite{tu2021variance} & 0.00 & 17.955(64.486)*  & 0.5833(0.3441)*  & 0.2267(0.0338) & 0.1560(0.0223)  \\
		& 0.10 & 5.5002(21.727)*  & 0.3885(0.0772)   & 0.2840(0.0507) & 0.1959(0.0287)  \\
		& 0.20 & 8.0043(53.507)*  & 0.5158(0.0862)   & 0.3720(0.0544) & 0.2713(0.0392) \\
		\hline \quad\quad\quad
		m=31	& $\alpha_n$ & $n=200$ & $n=300$ &   $n=500$ & $n=1000$ \\
		\hline \cite{2013Communication} & 0.00 & 0.2968(0.0401)   & 0.2069(0.0281)   & 0.1449(0.0201) & 0.0932(0.0128)  \\
		& 0.10 & 0.3731(0.0521)   & 0.3099(0.0378)   & 0.2713(0.0245) & 0.2432(0.0132)  \\
		& 0.20 & 0.5982(0.0615)   & 0.5306(0.0414)   & 0.4937(0.0265) & 0.4704(0.0149)  \\
		\cite{huo2019distributed} & 0.00 & 0.1908(0.0276)   & 0.1579(0.0225)   & 0.1234(0.0189) & 0.0866(0.0128)  \\
		& 0.10 & 0.4260(0.0720)   & 0.3477(0.0534)   & 0.2674(0.0432) & 0.1791(0.0197)  \\
		& 0.20 & 10.762(34.987) & 1.4455(0.3227)   & 0.9426(0.1766) & 0.6003(0.0886)  \\
		 ROSE(med) & 0.00 & 0.2923(0.0496)   & 0.2367(0.0369)   & 0.1807(0.0287) & 0.1250(0.0212)  \\
		& 0.10 & 0.3359(0.0561)   & 0.2636(0.0476)   & 0.2063(0.0347) & 0.1429(0.0207)  \\
		& 0.20 & 0.3867(0.0593)   & 0.3030(0.0523)   & 0.2331(0.0363) & 0.1645(0.0230)  \\
		& 0.00 & 0.2058(0.029)    & 0.1710(0.0245)   & 0.1305(0.0189) & 0.0911(0.0132)  \\
		& 0.10 & 0.2969(0.0455)   & 0.2314(0.0373)   & 0.1770(0.0256) & 0.1227(0.0177)  \\
		& 0.20 & 0.3501(0.0505)   & 0.2857(0.0447)   & 0.2151(0.0317) & 0.1518(0.0202)  \\
		\cite{Jordan2019Communication} & 0.00 & 104.78(590.15)*  & 0.6076(0.6602)*  & 0.1232(0.0193) & 0.0880(0.0137)  \\
		& 0.10 & 0.4536(0.0648)   & 0.3644(0.0583)   & 0.2757(0.0508) & 0.1967(0.0288)  \\
		& 0.20 & 8.1895(52.346)  & 1.4979(0.4350)   & 0.9897(0.2114) & 0.6281(0.1055)  \\
		\cite{tu2021variance} & 0.00 & 12.520(35.540)   & 0.8253(2.2017)   & 0.1310(0.0196) & 0.0927(0.0148)  \\
		& 0.10 & 6096.7(60865)*   & 0.2354(0.0418)   & 0.1726(0.0271) & 0.1232(0.0186)  \\
		& 0.20 & 0.4283(0.1549)   & 0.3218(0.0511)   & 0.2457(0.0345) & 0.1673(0.0246)  \\
			\hline \quad\quad\quad
		m=101  	& $\alpha_n$ & $n=200$ & $n=300$ &   $n=500$ & $n=1000$ \\
		\hline \cite{2013Communication} & 0.00 & 0.1991(0.0190)    & 0.1343(0.0145)   & 0.0867(0.0106) & 0.0547(0.0070)   \\
		& 0.10 & 0.2470(0.0287)   & 0.2315(0.0210)   & 0.2263(0.0106) & 0.2288(0.0071)  \\
		& 0.20 & 0.4897(0.0270)   & 0.4764(0.0214)   & 0.4689(0.0135) & 0.4643(0.0087)  \\
		\cite{huo2019distributed} & 0.00 & 0.1028(0.0156)   & 0.0879(0.0118)   & 0.0668(0.0102) & 0.0475(0.0069)  \\
		& 0.10 & 0.2367(0.0402)   & 0.1895(0.0332)   & 0.1461(0.0203) & 0.1046(0.0136)  \\
		& 0.20 & 0.8987(0.1675)   & 0.6787(0.1231)   & 0.4968(0.0676) & 0.3384(0.0525)  \\
		 ROSE(med) 	& 0.00 & 0.1647(0.0271)   & 0.1323(0.0212)   & 0.1026(0.0194) & 0.0717(0.0109)  \\
		& 0.10 & 0.1848(0.0313)   & 0.1490(0.0230)   & 0.1144(0.0156) & 0.0784(0.0106)  \\
		& 0.20 & 0.2162(0.0330)   & 0.1736(0.0278)   & 0.1371(0.0218) & 0.0956(0.0144)  \\
		 ROSE(K=10) & 0.00 & 0.1133(0.0165)   & 0.0933(0.0137)   & 0.0726(0.0111) & 0.0513(0.0071)  \\
		& 0.10 & 0.1678(0.0241)   & 0.1309(0.0185)   & 0.1014(0.0142) & 0.0702(0.0096)  \\
		& 0.20 & 0.2030(0.0298)   & 0.1669(0.0230)   & 0.1345(0.0175) & 0.0983(0.0106)  \\
		\cite{Jordan2019Communication} & 0.00 & 23.707(109.92)*  & 0.6154(0.9426)*  & 0.0701(0.0134) & 0.0468(0.0071)  \\
		& 0.10 & 0.2426(0.0422)   & 0.2001(0.0309)   & 0.1534(0.0221) & 0.1077(0.0149)  \\
		& 0.20 & 0.8853(0.1511)   & 0.7014(0.1283)   & 0.5184(0.0931) & 0.3464(0.0546)  \\
		\cite{tu2021variance} & 0.00 & 989.94(9476.7)*  & 0.4286(0.2333)*  & 0.0733(0.0116) & 0.0521(0.0077)  \\
		& 0.10 & 13.695(56.833)*  & 0.5328(4.0322)   & 0.0954(0.0162) & 0.0702(0.0115)  \\
		& 0.20 & 9.7351(62.207)*  & 0.1839(0.0244)   & 0.1418(0.0198) & 0.1006(0.0144)
		\\ \hline
	\end{tabular}
\begin{tablenotes}
	\footnotesize
	\item * means that the estimators do not converge, so we represent the estimation results after 3 rounds of iteration instead of the estimation results after 10 rounds.
\end{tablenotes}
\end{threeparttable}
}
\end{table}

\begin{table}
	\renewcommand\arraystretch{1.0}
	\centering
	\caption{Poisson regression model. }
	\label{T3}
	\scalebox{0.8}{
	\begin{threeparttable}
		\begin{tabular}{lccccc}
	\hline \quad\quad\quad m=11  	& $\alpha_n$ & $n=200$ & $n=300$ &   $n=500$ & $n=1000$ \\
			\hline \cite{2013Communication} & 0.00 & 0.1444(0.0217) & 0.1110(0.0195)
			& 0.0806(0.0124) & 0.0546(0.0080) \\
			& 0.10 & 0.1761(0.0191) & 0.1510(0.0149) & 0.1319(0.0089) & 0.1195(0.0054)\\
			& 0.20 & 0.2523(0.0153) & 0.2365(0.0125) & 0.2256(0.0075) & 0.2192(0.0048)\\
			\cite{huo2019distributed} & 0.00 & 0.1198(0.0192) & 0.0969(0.0162) & 0.0745(0.0110) & 0.0521(0.0078) \\
			& 0.10 & 0.1791(0.0287) & 0.1530(0.0211) & 0.1243(0.0155) & 0.0975(0.0118)\\
			& 0.20 & 0.2689(0.0392) & 0.2250(0.0280) & 0.1949(0.0231) & 0.1614(0.0137)\\
			 ROSE(med) & 0.00 & 0.2075(0.0393) & 0.1521(0.0281) & 0.1145(0.0187) & 0.0781(0,0121)  \\
			& 0.10 & 0.2201(0.0454) & 0.1713(0.0338) & 0.1250(0.0206) &0.0854(0.0138)  \\
			& 0.20 & 0.2700(0.0554) & 0.1952(0.0363) & 0.1398(0.0251) &0.0937(0.0175)  \\
			 ROSE(K=10) & 0.00 & 0.1373(0.0222) & 0.1065(0.0172) & 0.0827(0.0108) & 0.0561(0.0082) \\
			& 0.10 & 0.1537(0.0250) & 0.1227(0.0173) & 0.0910(0.0160) & 0.0621(0.0093) \\
			& 0.20 & 0.1897(0.0312) & 0.1426(0.0236) & 0.1048(0.0179) & 0.0693(0.0114) \\
			\cite{Jordan2019Communication} & 0.00 & 19.385(2408.5)* & 0.9744(2.8784)* & 0.3638(0.2712) & 0.0539(0.0079) \\
			& 0.10 & 3.7598(28.073)* & 1.2180(0.4729)* & 1.1698(0.0887) & 1.1018(0.0123)  \\
			& 0.20 & 12.074(1047.9)* & 1.1184(0.5134)* & 1.0531(0.0567) & 1.0217(0.0070)  \\
			\cite{tu2021variance} & 0.00 & 7.0880(44.641)* & 0.5319(0.3318)* & 0.2584(0.2261)&  0.0570(0.0096) \\
			& 0.10 & 3.3727(32.516)* & 1.2491(0.6251)* & 0.2666(0.2171) & 0.0635(0.0092) \\
			& 0.20 & 2.3816(1.7606)* & 1.2859(0.1606)* & 0.2057(0.1791) & 0.0662(0.0098) \\
			\hline \quad\quad\quad
			m=31	& $\alpha_n$ & $n=200$ & $n=300$ &   $n=500$ & $n=1000$ \\
			\hline \cite{2013Communication} & 0.00 & 0.0835(0.0102) & 0.0632(0.0084)
			& 0.0485(0.0072) & 0.0326(0.0500) \\
			& 0.10 & 0.1373(0.0099) & 0.1278(0.0086) & 0.1215(0.0052) &0.1173(0.0031)\\
			& 0.20 & 0.2371(0.0090) & 0.2321(0.0063) &  0.2303(0.0046) & 0.2274(0.0032)\\
			\cite{huo2019distributed} & 0.00 & 0.0689(0.0094) & 0.0567(0.0071) & 0.0435(0.0065) & 0.0314(0.0048) \\
			& 0.10 & 0.1215(0.0151) & 0.1056(0.0157) & 0.0890(0.0087) &0.0762(0.0050)\\
			& 0.20 & 0.1993(0.0228) & 0.1818(0.0178) &  0.1641(0.0134) & 0.1501(0.0086)\\
			 ROSE(med) & 0.00 & 0.1101(0.0191) & 0.0860(0.0126) & 0.0667(0.0115) &  0.0459(0.0067) \\
			& 0.10 & 0.1244(0.0202) & 0.0974(0.0140) &  0.0735(0.0119) & 0.0511(0.0079) \\
			& 0.20 & 0.1581(0.0033) & 0.1202(0.0209) &  0.0854(0.0134) & 0.0576(0.0086) \\
			 ROSE(K=10) & 0.00 & 0.0770(0.0116) & 0.0617(0.0087) & 0.0486(0.0068)& 0.0343(0.0053) \\
			& 0.10  & 0.0885(0.0135) & 0.0704(0.0102) & 0.0536(0.0079)& 0.0368(0.0055) \\
			& 0.20 & 0.1180(0.0207) & 0.0896(0.0144) &  0.0642(0.0117) & 0.0423(0.0074) \\
			\cite{Jordan2019Communication} & 0.00 & 5$\times10^4$(3$\times 10^{10}$)* & 0.7538(1.0732)* & 0.4450(0.2886) & 0.0343(0.0190) \\
			& 0.10 & 134.25(1314.3)* & 1.3128(0.2528)* & 1.1665(0.0902) & 1.1000(0.0031) \\
			& 0.20 & 1$\times 10^4$(1$\times 10^9$)* & 1.0820(0.1929)* & 1.0529(0.0605) & 1.0182(0.0043) \\
			\cite{tu2021variance} & 0.00 & 2.0005(1.8288)* & 0.5583(0.3016)*& 0.3417(0.2647) & 0.0345(0.0088)  \\
			& 0.10 & 2.8344(4.5085)* & 1.3007(0.2011)* & 0.3377(0.2322) & 0.0375(0.0051) \\
			& 0.20 & 6.7421(40.097)* & 1.2596(0.2051)* & 0.3174(0.2636) & 0.0436(0.0058) \\
			\hline \quad\quad\quad
			m=101  	& $\alpha_n$ & $n=200$ & $n=300$ &   $n=500$ & $n=1000$ \\
			\hline \cite{2013Communication} & 0.00 & 0.0476(0.0070) & 0.0359(0.0052)
			& 0.0268(0.0037) & 0.0183(0.0025) \\
			& 0.10 & 0.1191(0.0048) & 0.1177(0.0037) & 0.1176(0.0025)  & 0.1162(0.0018)\\
			& 0.20 & 0.2309(0.0048) & 0.2324(0.0033) & 0.2311(0.0023)  & 0.2308(0.0015)\\
			\cite{huo2019distributed} & 0.00 & 0.0387(0.0063) & 0.0313(0.0047) & 0.0248(0.0036) & 0.0176(0.0026) \\
			& 0.10 & 0.0868(0.0072) & 0.0778(0.0057) & 0.0728(0.0040)  &0.0668(0.0029)\\
			& 0.20 & 0.1595(0.0117) & 0.1564(0.0089) & 0.1474(0.0064)  &0.1433(0.0041)\\
			 ROSE(med) & 0.00 &0.0592(0.0105) & 0.0469(0.0076) & 0.0370(0.0052) & 0.0250(0.0039)\\
			& 0.10 &  0.0678(0.0096)& 0.0545(0.0088) & 0.0400(0.0061) &0.0288(0.0049) \\
			& 0.20 & 0.0865(0.0125) & 0.0685(0.0103) & 0.0493(0.0085) &0.0322(0.0051) \\
			 ROSE(K=10) & 0.00 & 0.0419(0.0067) & 0.0332(0.0046) & 0.0264(0.0039) &0.0181(0.0029) \\
			& 0.10 & 0.0498(0.0072) & 0.0404(0.0063) & 0.0300(0.0045) &0.0213(0.0034)  \\
			& 0.20 & 0.0704(0.0127) & 0.0555(0.0111) & 0.0386(0.0053) &0.0254(0.0037)  \\
			\cite{Jordan2019Communication} & 0.00 & 85.734(5$\times10^4$)* & 0.6631(0.4454)* & 0.4458(0.2948) & 0.0184(0.0058) \\
			& 0.10 & 6.9107(213.98)* & 1.2808(0.9510)* & 1.1843(0.0919) & 1.0988(0.0022) \\
			& 0.20 & 1$\times 10^{7}$(1$\times 10^{15}$)* & 1.1175(0.6231)* & 1.0493(0.0591) & 1.0167(0.0036) \\
			\cite{tu2021variance} & 0.00 & 4.0250(9.8595)*& 0.5451(0.3556)*&
			0.3735(0.2622) &0.0195(0.0049)
			\\
			& 0.10 & 2$\times10^4$(2$\times10^5$)* & 1.3810(0.3699)* & 0.3817(0.2573)  & 0.0270(0.0047)
			\\
			& 0.20 & 4.1269(17.288)* & 1.2940(0.2071)* & 0.3284(0.2736)  & 0.0256(0.0034)
			\\ \hline
		\end{tabular}
	\begin{tablenotes}
		\footnotesize
		\item * means that the estimators do not converge, so we represent the estimation results after 3 rounds of iteration instead of the estimation results after 10 rounds.
	\end{tablenotes}
\end{threeparttable}
	}
\end{table}
\subsubsection{Further comparison: when the samples on the central processor are anomalous}\label{5.1.3}
 We mainly compare three different parameter estimators: (1)  ROSE (RV, K=10), computed by Algorithm \ref{alg2}, which uses $\boldsymbol{\hat{\theta}}{'}_{vr}$, $\widehat{\nabla F}{'}_{vr}(\boldsymbol{\hat{\theta}}{'}_{vr})$ and $\widehat{\nabla^2 F}{'}_{vr}(\boldsymbol{\hat{\theta}}{'}_{vr})$ as the initial value, the gradient and the Hessian matrix respectively, and then compute $\boldsymbol{\hat{\theta}}{'}_{ros}$
with $K=10$; (2)  ROSE (K=10), calculated by Algorithm \ref{alg1} with $K=10$; and (3) the robust approximation based on the VRMOL estimator with $K=10$ (proposed by \cite{tu2021variance}). For logical model,  in the setting of  Byzantine machines, we  consider two cases: (1) the data are generated from the model (\ref{log2}), while for Byzantine machines, the response variable $Y$ in the original data is recorded as $1-Y$, and  Byzantine machines use the incorrect $Y$ in the following computations; and  (2) when we distribute normal dataset on machines, the $\boldsymbol{X}$ received by the Byzantine machine is $10$ times larger than the correct $\boldsymbol{X}$, and the central processor is also a Byzantine machine.
 And for poission model, we  consider: (1) the same as the simulations in Section \ref{5.1.2}, but the central processor is a Byzantine machine; and  (2) the same as case (2) for the logical regression model. It is worth noting that in case (1), the  Hessian matrix estimator is slight influenced,  but in case (2), the estimator is greatly influenced.

Tables \ref{T4}--\ref{T7} showcase simulation outcomes for logical  and Poisson regression models under the two Byzantine machine settings. All other settings except for Byzantine machines are the same as those in the previous simulations. 

In Table \ref{T4}, we observe that the averages of  root  squared errors resulting from  ROSE (RV,K=10) and  ROSE (K=10) are nearly identical. This is because the Hessian matrix estimation is only slightly affected by Byzantine machines. Conversely, the estimators provided by \cite{tu2021variance} fail to converge, resulting in their exclusion from the table. In Table \ref{T6}, since the Hessian  matrix estimator is  slightly influenced,  the results obtained by the three methods are not significantly different from the corresponding results in Table \ref{T3}. However, it is evident that among these three methods,  ROSE (RV,K=10) provides the best outcome.

As for Tables \ref{T5} and \ref{T7}, since the Hessian matrice exhibit  significant impacts,  ROSE (K=10) performs very poor and we omit the corresponding results. Additionally, the estimators provided by \cite{tu2021variance} produce much worse outcomes than  ROSE (RV,K=10).

%
\begin{table}
	\renewcommand\arraystretch{1.0}
	\centering
	\caption{Logistic regression model: case (1) of the Byzantine setting.}
	\label{T4}
	\scalebox{0.8}{
		\begin{tabular}{lccccc}
		\hline \quad\quad\quad m=11  	& $\alpha_n$ & $n=200$ & $n=300$ &   $n=500$ & $n=1000$ \\
			\hline
			 ROSE(RV,K=10)
			& 0.10 & 0.3549(0.0558) & 0.3102(0.0408) & 0.2388(0.0454) &0.1641(0.0243)  \\
			& 0.20 & 0.3890(0.0591) & 0.3158(0.0473) & 0.2532(0.0305) &0.1792(0.0263) \\
			 ROSE(K=10)
			& 0.10 & 0.3389(0.0496) & 0.2956(0.0369) & 0.2317(0.0416) & 0.1591(0.0227) \\
			& 0.20 & 0.3748(0.0550) & 0.3088(0.0428) & 0.2475(0.0271) & 0.1745(0.0246) \\
			\hline \quad\quad\quad
			m=31	& $\alpha_n$ & $n=200$ & $n=300$ &   $n=500$ & $n=1000$ \\
			\hline
			 ROSE(RV,K=10)
			& 0.10 & 0.2177(0.0302) & 0.1797(0.0269) & 0.1409(0.0224) &0.1020(0.0139)  \\
			& 0.20 & 0.2524(0.0366) & 0.2066(0.0236) & 0.1656(0.0198) &0.1204(0.0143) \\
			 ROSE(K=10)
			& 0.10 & 0.2129(0.0292) & 0.1759(0.0231) & 0.1378(0.0204) & 0.1008(0.0124) \\
			& 0.20 & 0.2493(0.0340) & 0.2070(0.0217) & 0.1669(0.0168) & 0.1211(0.0127) \\
			\hline \quad\quad\quad
			m=101  	& $\alpha_n$ & $n=200$ & $n=300$ &   $n=500$ & $n=1000$ \\
			\hline
			 ROSE(RV,K=10)
			& 0.10 & 0.1320(0.0179) & 0.1081(0.0129) & 0.0862(0.0099) &0.0607(0.0074)  \\
			& 0.20 & 0.1752(0.0159)&0.1456(0.0118)
			&0.1179(0.0107) &0.0860(0.0076) \\
			 ROSE(K=10)
			& 0.10 & 0.1341(0.0158) & 0.1095(0.0116) & 0.0885(0.0097) & 0.0617(0.0066) \\
			& 0.20 & 0.1847(0.0141) & 0.1531(0.0108) & 0.1248(0.0092) & 0.0922(0.0070)
			\\ \hline
		\end{tabular}
	}
\end{table}
\begin{table}
	\renewcommand\arraystretch{1.0}
	\centering
	\caption{Logistic regression model: case (2) of the Byzantine setting. }
	\label{T5}
	\scalebox{0.8}{
		\begin{tabular}{lccccc}
			\hline \quad\quad\quad m=11  	& $\alpha_n$ & $n=200$ & $n=300$ &   $n=500$ & $n=1000$ \\
			\hline
			 ROSE(RV,K=10)
			& 0.10 & 0.3609(0.0575) & 0.3031(0.0443) & 0.2251(0.0349) &0.1612(0.0223)  \\
			& 0.20 & 0.4030(0.0592) & 0.3320(0.0595) & 0.2549(0.0374) &0.1770(0.0242) \\
			\cite{tu2021variance}
			& 0.10 & 0.4754(0.0237) & 0.4611(0.0190) & 0.4578(0.1510) & 0.4485(0.0131) \\
			& 0.20 & 0.4882(0.0153) & 0.4878(0.0109) & 0.4837(0.0097) & 0.4759(0.0077) \\
			\hline \quad\quad\quad
			m=31	& $\alpha_n$ & $n=200$ & $n=300$ &   $n=500$ & $n=1000$ \\
			\hline
			 ROSE(RV,K=10)
			& 0.10 & 0.2185(0.0330) & 0.1825(0.0260) & 0.1429(0.0203) &0.1002(0.0154)  \\
			& 0.20 & 0.2533(0.0321) & 0.2052(0.0233) & 0.1626(0.0204) &0.1181(0.0140) \\
			\cite{tu2021variance}
			& 0.10 & 0.4548(0.0197) & 0.4577(0.0173) & 0.4507(0.0141) & 0.4462(0.0115) \\
			& 0.20 & 0.4890(0.0088) & 0.4856(0.0086) & 0.4833(0.0067) & 0.4772(0.0053) \\
			\hline \quad\quad\quad
			m=101  	& $\alpha_n$ & $n=200$ & $n=300$ &   $n=500$ & $n=1000$ \\
			\hline
			 ROSE(RV,K=10)
			& 0.10 & 0.1300(0.0186) & 0.1083(0.0145) & 0.0824(0.0103) &0.0587(0.0082)  \\
			& 0.20 & 0.1717(0.0147)&0.1416(0.0123)
			&0.1131(0.0088) &0.0827(0.0074) \\
			\cite{tu2021variance}
			& 0.10 & 0.4513(0.0165)& 0.4499(0.0138) & 0.4499(0.0130) & 0.4422(0.0106) \\
			& 0.20 & 0.4869(0.0082) & 0.4855(0.0057) & 0.4833(0.0052) & 0.4771(0.0051)
			\\ \hline
		\end{tabular}
	}
\end{table}
\begin{table}
	\renewcommand\arraystretch{1.0}
	\centering
	\caption{Poisson regression model: case (1) of the Byzantine setting.}
	\label{T6}
	\scalebox{0.8}{
	\begin{threeparttable}
		\begin{tabular}{lccccc}
			\hline \quad\quad\quad m=11  	& $\alpha_n$ & $n=200$ & $n=300$ &   $n=500$ & $n=1000$ \\
			\hline
			 ROSE(RV,K=10)
			& 0.10 & 0.1462(0.0230) & 0.1138(0.0158) & 0.0885(0.0144) &0.0606(0.0107)  \\
			& 0.20 & 0.1749(0.0288) & 0.1359(0.0205) & 0.1021(0.0172) &0.0687(0.0109) \\
			 ROSE(K=10)
			& 0.10 & 0.1603(0.0311) & 0.1258(0.0204) & 0.0942(0.0157) & 0.0650(0.0112) \\
			& 0.20 & 0.2071(0.0541) & 0.1619(0.0268) & 0.1173(0.0192) & 0.1072(0.0203) \\
			\cite{tu2021variance}
			& 0.10 & 16.472(32.739)* & 1.7764(1.2423)* & 0.4493(0.2860) & 0.0939(0.0624) \\
			& 0.20 & 15.634(60.638)* & 1.5296(0.9660)* & 0.3805(0.2289) & 0.0905(0.0370) \\
			\hline \quad\quad\quad
			m=31	& $\alpha_n$ & $n=200$ & $n=300$ &   $n=500$ & $n=1000$ \\
			\hline
		 ROSE(RV,K=10)
		& 0.10 & 0.0873(0.0126) & 0.0687(0.0102) & 0.0521(0.0084) &0.0366(0.0052)  \\
		& 0.20 & 0.1059(0.0154) & 0.0838(0.0126) & 0.0612(0.0089) &0.0435(0.0063) \\
		 ROSE(K=10)
		& 0.10 & 0.1025(0.0212) & 0.0800(0.0229) & 0.0595(0.0109) & 0.0411(0.0062) \\
		& 0.20 & 0.1439(0.0497) & 0.1102(0.0255) & 0.0765(0.0124) & 0.0490(0.0087) \\
		\cite{tu2021variance}
		& 0.10 & 26.991(110.99)* & 1.7151(1.5210)* & 0.3488(0.2285) & 0.0619(0.0397) \\
		& 0.20 & 11.272(38.896)* & 1.5116(0.7891)* & 0.3322(0.2220) & 0.0648(0.0264) \\
			\hline \quad\quad\quad
			m=101  	& $\alpha_n$ & $n=200$ & $n=300$ &   $n=500$ & $n=1000$ \\
			\hline
			 ROSE(RV,K=10)
		& 0.10 & 0.0485(0.0078) & 0.0386(0.0054) & 0.0295(0.0047) &0.0199(0.0032)  \\
		& 0.20 & 0.0610(0.0080)&0.0487(0.0072)
		&0.0367(0.0050) &0.0253(0.0337) \\
		 ROSE(K=10)
		& 0.10 & 0.0620(0.0153) & 0.0482(0.0092) & 0.0376(0.0080) & 0.0251(0.0044) \\
		& 0.20 & 0.0974(0.0314) & 0.0779(0.0216) & 0.0540(0.0118) & 0.0339(0.0107) \\
		\cite{tu2021variance}
		& 0.10 & 9.1160(15.198)* & 1.3931(0.7817)* & 0.3500(0.2110) & 0.0460(0.0395) \\
		& 0.20 & 7.4089(12.182)* & 1.4865(0.7846)* & 0.2926(0.2238) & 0.0586(0.0478)
			\\ \hline
		\end{tabular}
	\begin{tablenotes}
		\footnotesize
		\item * means that the estimators do not converge, so we represent the estimation results after 2 rounds of iteration instead of the estimation results after 10 rounds.
	\end{tablenotes}
\end{threeparttable}
	}
\end{table}
\begin{table}
	\renewcommand\arraystretch{1.0}
	\centering
	\caption{Poisson regression model: case (2) of the Byzantine setting.}
	\label{T7}
	\scalebox{0.8}{
		\begin{tabular}{lccccc}
			\hline \quad\quad\quad m=11  	& $\alpha_n$ & $n=200$ & $n=300$ &   $n=500$ & $n=1000$ \\
			\hline
			 ROSE(RV,K=10)
			& 0.10 & 0.1578(0.0236) & 0.1243(0.0204) & 0.0974(0.0151) &0.0656(0.0095)  \\
			& 0.20 & 0.2217(0.0425) & 0.1719(0.0296) & 0.1232(0.0202) &0.0850(0.0123) \\
			\cite{tu2021variance}
			& 0.10 & 0.4988(0.0070) & 0.5007(0.0046) & 0.5023(0.0040) & 0.5041(0.0027) \\
			& 0.20 & 0.5044(0.0061) & 0.5057(0.0048) & 0.5065(0.0036) & 0.5068(0.0026) \\
			\hline \quad\quad\quad
			m=31	& $\alpha_n$ & $n=200$ & $n=300$ &   $n=500$ & $n=1000$ \\
			\hline
			 ROSE(RV,K=10)
			& 0.10 & 0.0972(0.0138) & 0.0759(0.0107) & 0.0580(0.0093) &0.0410(0.0057)  \\
			& 0.20 & 0.1401(0.0238) & 0.1087(0.0166) & 0.0815(0.0122) &0.0546(0.0089) \\
			\cite{tu2021variance}
			& 0.10 & 0.4983(0.0066) & 0.5006(0.0047) & 0.5025(0.0038) & 0.5034(0.0025) \\
			& 0.20 & 0.5047(0.0050) & 0.5057(0.0046) & 0.5064(0.0035) & 0.5065(0.0023) \\
			\hline \quad\quad\quad
			m=101  	& $\alpha_n$ & $n=200$ & $n=300$ &   $n=500$ & $n=1000$ \\
			\hline
			 ROSE(RV,K=10)
			& 0.10 & 0.0541(0.0077) & 0.0435(0.0062) & 0.0333(0.0055) &0.0227(0.0035)  \\
			& 0.20 & 0.0821(0.0133)&0.0645(0.0095)
			&0.0475(0.0072) &0.0302(0.0042) \\
			\cite{tu2021variance}
			& 0.10 & 0.4980(0.0054)& 0.5013(0.0047) & 0.5016(0.0037) & 0.5040(0.0027) \\
			& 0.20 & 0.5045(0.0055) & 0.5056(0.0037) & 0.5062(0.0029) & 0.5069(0.0025)
			\\ \hline
		\end{tabular}
	}
\end{table}

\subsection{ The real data example}
The dataset is a subset of the NIST (National Institute of Standards and Technology) dataset, which can be accessed at \url{http://yann.lecun.com/exdb/mnist/}. It consists of  training and testing sets, with the former containing 60,000 images and labels, and the latter including 10,000 images and labels. Each image is a gray handwritten digital image with $28\times28$ pixels ranging from 0 to 9, with white characters on a black background and pixel values ranging from 0 to 255.

We select three digits, 6, 8, and 9, which are known to be more challenging to distinguish. Subsequently, we aim to  train two logistic classifiers for $6,9$ and $8,9$, respectively. Initially, we eliminate variables that have 75\% of observations equaling zero and then employ the Lasso-logistic regression introduced in \cite{2010Regularization} to filter the remaining variables. Ultimately, we select 15 and 26 key variables from  a total of 784 variables for the two classifiers respectively. We  use the classification error on the testing set to evaluate the performance of different parameter estimation methods. The original data set was used for normal machines, while for Byzantine machines, we replace the independent variables corresponding to digits 6, 8, and 9 with the independent variables corresponding to digits 0 and 1. We consider two proportions of Byzantine machines, $\alpha_n= 0.00$ and $0.20$. The results are assembled in Tables \ref{T8} and \ref{T9}.

From Tables \ref{T8} and \ref{T9}, we can see that when there is no Byzantine machine, the results of all the methods are comparable, and  when 20\% machines are Byzantine machines, the prediction accuracies of the methods proposed by \cite{2013Communication} and \cite{huo2019distributed} are significantly lower than the other competitors. In the classification of 8 and 9, when $m=30$ and $n=396$,  ROSE (med) and  ROSE (K=5) perform much better than \cite{tu2021variance}. Moreover, the method proposed by \cite{Jordan2019Communication} does not converge when there exist 20\% Byzantine machines because the gradient value cannot converge to $0$. 

\begin{table}
	\renewcommand\arraystretch{1.0}
	\centering
	\caption{Prediction accuracy.}
	\label{T8}
	\scalebox{1.0}{
		\begin{tabular}{lccccc}
				\hline Classification of 6 and 9  & & $m=10$ & $m=20$ & $m=30$  \\
			\quad\quad(15 variables)	& $\alpha_n$ & $n=1188$ & $n=594$ &   $n=396$  \\
			\hline \cite{2013Communication} & 0.00 & 90.95\% & 90.80\%
			& 90.85\%  \\
			& 0.20 & 89.88\% & 89.83\% & 89.88\% \\
			\cite{huo2019distributed} & 0.00 & 90.90\% & 90.90\% & 90.85\%  \\
			& 0.20 & 51.45\% & 44.13\% & 36.20\%  \\
			 ROSE(med) & 0.00 & 90.80\% & 90.70\% & 90.59\%  \\
			& 0.20 & 90.95\% & 91.00\% & 91.05\%  \\
			 ROSE(K=5) & 0.00 & 90.80\% & 90.90\% & 90.75\% \\
			& 0.20 & 90.95\% & 91.10\% & 91.10\%  \\
			\cite{Jordan2019Communication} & 0.00 & 90.90\% & 90.85\% & 90.85\% \\
			& 0.20 & --- & --- & ---  \\
			\cite{tu2021variance} & 0.00 & 90.85\% & 90.85\% & 90.80\%  \\
			& 0.20 & 91.05\% & 91.00\% & 90.75\% \\
			\hline
		\end{tabular}
	}
\end{table}

\begin{table}
	\renewcommand\arraystretch{1.0}
	\centering
	\caption{Prediction accuracy.}
	\label{T9}
	\scalebox{1.0}{
		\begin{tabular}{lccccc}
				\hline Classification of 8 and 9  & & $m=10$ & $m=20$ & $m=30$  \\
			\quad\quad(26 variables)	& $\alpha_n$ & $n=1188$ & $n=594$ &   $n=396$  \\
			\hline \cite{2013Communication} & 0.00 & 87.24\% & 87.34\%
			& 87.09\%  \\
			& 0.20 & 78.87\% & 76.20\% & 74.38\% \\
			\cite{huo2019distributed} & 0.00 & 87.14\% & 87.14\% & 87.14\%  \\
			& 0.20 & 82.96\% & 79.48\% & 73.17\%  \\
			 ROSE(med) & 0.00 & 87.49\% & 87.14\% & 86.94\%  \\
			& 0.20 & 87.04\% & 86.89\% & 86.28\%  \\
			 ROSE(K=5) & 0.00 & 87.44\% & 87.19\% & 87.29\% \\
			& 0.20 & 87.34\% & 87.19\% & 87.04\%  \\
			\cite{Jordan2019Communication} & 0.00 & 87.14\% & 87.34\% & 83.01\% \\
			& 0.20 & --- & --- & ---  \\
			\cite{tu2021variance} & 0.00 & 87.49\% & 87.04\% & 80.53\%  \\
			& 0.20 & 87.24\% & 87.14\% & 63.54\% \\
			\hline
		\end{tabular}
	}
\end{table}
\section{Discussions}\label{sec6}
We have established the asymptotic normality and convergence rate for the proposed estimator, which also exhibits robustness when Byzantine machines comprise a relatively minor portion of the distributed network. To achieve the robustness of the estimator,  other methods such as trimmed mean (\cite{2018Byzantine}) and geometric median (\cite{2014Distributed}), can be also utilized to estimate initial parameter values for Newton-Raphson iteration. Besides, if minimizing communication costs is paramount, the local estimator of the central processor can directly serve as an initial parameter value. Naturally, the convergence rate of the final one-step estimator may  depend on  initial estimator.

 Some research directions  remain open and interesting. First, our method mainly applies to M-estimation problems  that satisfy loss second-order differentiability, rendering these techniques non-applicable when dealing with quantile regression models devoid of empirical Hessian matrices.  Besides,
 it cannot solve the problem of non-parametric distributed estimation. Second, we do not consider instances where  $p$ diverges at a faster rate, such as $p/n\to c>0$ or $p>n$ but $p^3=o(mn)$. Additionally, we only provide the convergence rate when $m=o(n)$, suggesting that when $m$ diverges faster than $n$, we may be
 able to improve the estimator to have another convergence rate. Third, this method could pave the way for enhancing estimation robustness while negligibly impacting estimation accuracy in decentralized distributed problems or streaming datasets. Last,   ROSE can be applied to model checking. Due to the asymptotic normality of this estimator requiring a relatively small proportion of Byzantine machines, it can be used to judge whether the data blocks that do not obey a certain model setting in multiple data blocks exceed a certain proportion.

\section{Assumptions}\label{sec7}
\begin{assumption}(Parameter space)\label{a1}
	The parameter space $\boldsymbol{\Theta} \subset \mathbb{R}^p$ is a compact convex set, and $\boldsymbol{\theta^*}$ is an interior point in $\boldsymbol{\Theta}$. The $\ell_2$-radius $D=\max_{\boldsymbol{\theta}\in \boldsymbol{\Theta}}\left\|\boldsymbol{\theta}-\boldsymbol{\theta^*}\right\|_2$ is bounded.
\end{assumption}
\begin{assumption}\label{a2}(Convexity)
	The loss function $f(\boldsymbol{x},\boldsymbol{\theta})$ is convex with respect to $\boldsymbol{\theta}\in \boldsymbol{\Theta}$ for all $\boldsymbol{x}$.
\end{assumption}
\begin{assumption}\label{a3}
	(Bounded eigenvalue for Hessian matrix)	The loss function $f(\boldsymbol{x},\boldsymbol{\theta})$ is twice differentiable, and there exists a positive constant $\lambda$ such that $\lambda^{-1}\leq \lambda_{\min}(\nabla^2 F_{\mu}(\boldsymbol{\theta^*}))\leq \lambda_{\max}(\nabla^2 F_{\mu}(\boldsymbol{\theta^*}))\leq \lambda$.
\end{assumption}
\begin{assumption}\label{a4}(Lipschitz continuous)
	There exists a function $L(\boldsymbol{X})$ and a positive constant $\delta$ such that for arbitrary $\boldsymbol{\theta}_1, \boldsymbol{\theta}_2\in B(\boldsymbol{\theta^*}, \delta)\subset \boldsymbol{\Theta}$,
	$\left\|\nabla^2f(\boldsymbol{X},\boldsymbol{\theta}_1)-\nabla^2f(\boldsymbol{X},\boldsymbol{\theta}_2)\right\|\leq  \sqrt{p}L(\boldsymbol{X})\left\|\boldsymbol{\theta}_1-\boldsymbol{\theta}_2\right\|$, where $L(\boldsymbol{X})$ satisfies $\mathbb{E}[L^4(\boldsymbol{X})]\leq L^4$ for some constant $L>0$.
\end{assumption}
\begin{assumption}\label{a5}(Sub-exponential for the gradient)
	There exists a positive constant $t$  such that
$$		\max\limits_{1\leq l\leq p}\mathbb{E}[\operatorname{exp}\{t|\nabla_{\theta_l}f(\boldsymbol{X},\boldsymbol{\theta^*})|\}]\leq 2.$$
\end{assumption}
\begin{assumption}\label{a6}(Bounded eigenvalue)
	There exists a positive constant $\lambda_0$ such that
	$$\lambda_0^{-1}\leq\|\mathbb{E}[\{\nabla f(\boldsymbol{X},\boldsymbol{\theta^*})\}^{\otimes2}]\|\leq \lambda_0.$$
\end{assumption}
\begin{assumption}\label{a7}(Sub-exponential for entries in the Hessian matrix)
	There exists a positive constant $t$  such that 
	$$\max\limits_{1\leq l_1,l_2\leq p}\mathbb{E}[\operatorname{exp}\{t|\nabla^2_{\theta_{l_1}\theta_{l_2}}f(\boldsymbol{X},\boldsymbol{\theta^*})-\nabla^2_{\theta_{l_1}\theta_{l_2}}F_{\mu}(\boldsymbol{\theta^*})|\}]\leq 2.$$
\end{assumption}
\begin{assumption}\label{a8}
	(Sub-exponential for the inner product of vectors in the Hessian matrix)
	Let $[\nabla^2 F_{\mu}(\boldsymbol{\theta^*})]^{-1}_{l\cdot}$ be the $l$-th row of $[\nabla^2 F_{\mu}(\boldsymbol{\theta^*})]^{-1}$ and $\nabla^2 f(\boldsymbol{X},\boldsymbol{\theta})_{\cdot l}$ be the $l$-th column of $\nabla^2 f(\boldsymbol{X},\boldsymbol{\theta})$.
	There exists two positive constants $t$ and $\delta$ such that for any $l_1,l_2\in[p]$, if $\boldsymbol{\theta}\in B(\boldsymbol{\theta^*},\delta)$, $$\mathbb{E}[\exp(t|\langle[\nabla^2 F_{\mu}(\boldsymbol{\theta^*})]^{-1}_{l_1\cdot},\nabla^2 f(\boldsymbol{X},\boldsymbol{\theta})_{\cdot l_2}\rangle|)]\leq 2.$$
\end{assumption}
\begin{assumption}\label{a9}
	(Moment ratio restriction 1)\label{berry1} For $l\in[p]$, there exists a positive constant $R_b$ such that
	\begin{align*}
		\frac{\mathbb{E}[|\langle[\nabla^2 F_{\mu}(\boldsymbol{\theta^*})]^{-1}_{l\cdot},\nabla f(\boldsymbol{X},\boldsymbol{\theta^*})\rangle|^3]}{\mathbb{E}[\langle[\nabla^2 F_{\mu}(\boldsymbol{\theta^*})]^{-1}_{l\cdot},\nabla f(\boldsymbol{X},\boldsymbol{\theta^*})\rangle^2]}\leq R_b.
	\end{align*}
\end{assumption}
\begin{assumption}\label{a10}(Moment ratio restriction 2)
	For any $l\in[p]$, there exist two positive constants $R_g$ and $\delta$ such that for any $\boldsymbol{\theta}\in B(\boldsymbol{\theta^*},\delta)$,
	\begin{align*}
		\frac{\mathbb{E}[|\nabla_{\theta_l} f(\boldsymbol{X},\boldsymbol{\theta})-\nabla_{\theta_l} F_{\mu}(\boldsymbol{\theta})|^3]}{\mathbb{E}[\{\nabla_{\theta_l} f(\boldsymbol{X},\boldsymbol{\theta})-\nabla_{\theta_l} F_{\mu}(\boldsymbol{\theta})\}^2]}\leq R_g.
	\end{align*}
\end{assumption}
\begin{assumption}\label{a102}(Moment ratio restriction 3)
	For any $l_1,l_2\in[p]$, there exist two positive constants $R_h$ and $\delta$ such that for any $\boldsymbol{\theta}\in B(\boldsymbol{\theta^*},\delta)$
	\begin{align*}
		\frac{\mathbb{E}[|\nabla^2_{\theta_{l_1}\theta_{l_2}} f(\boldsymbol{X},\boldsymbol{\theta})-\nabla^2_{\theta_{l_1}\theta_{l_2}} F_{\mu}(\boldsymbol{\theta})|^3]}{\mathbb{E}[\{\nabla^2_{\theta_{l_1}\theta_{l_2}} f(\boldsymbol{X},\boldsymbol{\theta})-\nabla^2_{\theta_{l_1}\theta_{l_2}} F_{\mu}(\boldsymbol{\theta})\}^2]}\leq R_h.
	\end{align*}
\end{assumption}
\begin{assumption}(Sub-exponential for each coordinate of the gradient)\label{a11}
	There exist two positive constants $t$ and $\delta$ such that for any $\boldsymbol{\theta}\in B(\boldsymbol{\theta^*},\delta)$,
	\begin{align*}
		\mathbb{E}\left[\exp\left\{\frac{t|\nabla_{\theta_l} f(\boldsymbol{X},\boldsymbol{\theta})-\nabla_{\theta_l} f(\boldsymbol{X},\boldsymbol{\theta^*})-\nabla_{\theta_l} F_{\mu}(\boldsymbol{\theta})+\nabla_{\theta_l} F_{\mu}(\boldsymbol{\theta^*})|}{\|\boldsymbol{\theta}-\boldsymbol{\theta^*}\|}\right\}\right]\leq 2.
	\end{align*}
\end{assumption}
\begin{assumption}\label{a12}(Smooth)
	There exist two positive constants $C_H$ and $\delta$ such that for arbitrary $\boldsymbol{\theta}_1, \boldsymbol{\theta}_2\in B(\boldsymbol{\theta^*}, \delta)\subset \boldsymbol{\Theta}$,
	$\left\|\nabla^2 F_{\mu}(\boldsymbol{\theta}_1)-\nabla^2 F_{\mu}(\boldsymbol{\theta}_2)\right\|\leq  C_H\left\|\boldsymbol{\theta}_1-\boldsymbol{\theta}_2\right\|$.
\end{assumption}
\begin{remark}
	Assumptions \ref{a1} and  \ref{a2}
	are standard in classical statistical analysis of M-estimator(e.g. \cite{2000Asymptotic}). Assumption \ref{a3} ensures the strong local convexity of the loss function, and we can find similar assumptions in \cite{2013Communication}, \cite{Jordan2019Communication} and \cite{tu2021variance}. Assumption \ref{a4} guarantees the smoothness of the second-order derivatives of the loss function. {Similar assumptions can be found in \cite{Jordan2019Communication} and \cite{huo2019distributed}, however,  since we allow $p$ to diverge to infinity,
	 we add a factor $\sqrt{p}$  to $\|\boldsymbol{\theta}_1-\boldsymbol{\theta}_2\|$.}
	 Assumption \ref{a5} requires each entry of the gradient to follow a sub-exponential distribution, which is equivalent to Assumption F in \cite{tu2021variance}  and  weaker than  the sub-gaussian assumption in
	 \cite{2019Defending}. Assumption \ref{a6} requires that the covariance matrix of the gradient at the true value of the parameter has bounded eigenvalues, which is a regularity condition in research about M-estimators. Assumption \ref{a7} requires that each entry of the second derivative of the loss function at the true value of the parameter obeys a sub-exponential distribution, similar to Assumption  \ref{a5}. Assumption \ref{a8} requires that the inner product of the row vector of the inverse of the Hessian matrix and the second partial derivative of the loss function obey a sub-exponential distribution. Because, in many cases, the non-diagonal entries of the inverse of a Hessian matrix are small, this assumption is not difficult to satisfy. Assumptions \ref{a9}-\ref{a102} require the third-order moment of some random variables to be controlled by a constant multiple of their second-order moment, which is often adopted when using Berry-Esseen theorem to prove asymptotic normality.
	 It is worth noting that if a random variable has pseudo-independence or follows an elliptical distribution(e.g., \cite{Cui2018TEST}), then its third-order moment can be controlled by a constant multiple of its second-order moment. Assumptions \ref{a11} and \ref{a12} are two smoothness assumptions, which are similar to Assumptions C and D in \cite{tu2021variance}. For all assumptions, we do not directly assume the dimension $p$ of parameter $\boldsymbol{\theta}$ or the dimension $q$ of $\boldsymbol{X}$, and the requirements for $p$ and $q$ are hidden in the functions related to the loss function and its derivative.
\end{remark}
\begin{remark}
	Under Assumptions \ref{a3} and \ref{a6}, we can conclude that for any $l\in[p]$,
	\begin{align*}
		(\lambda^2\lambda_0)^{-1}\leq\sigma_l^2(\boldsymbol{\theta^*})\leq \lambda^2\lambda_0.
	\end{align*}
\end{remark}
\bibliographystyle{apalike}

\bibliography{reference(d)}

\begin{thebibliography}{}

\bibitem[Battey et~al., 2018]{2018Distributed}
Battey, H., Fan, J., Liu, H., Lu, J., and Zhu, Z. (2018).
\newblock Distributed testing and estimation under sparse high dimensional
  models.
\newblock {\em The Annals of Statistics}, 46:1352--1382.

\bibitem[Cui et~al., 2018]{Cui2018TEST}
Cui, H., Guo, W., and Zhong, W. (2018).
\newblock Test for high-dimensional regression coefficients using refitted
  cross-validation variance estimation.
\newblock {\em The Annals of Statistics}, 46:958--988.

\bibitem[Duan et~al., 2020]{2020Learning}
Duan, R., Luo, C., Schuemie, M.~J., Tong, J., Liang, C.~J., Chang, H.~H.,
  Regina, B.~M., Bian, J., Xu, H., and Holmes, J.~H. (2020).
\newblock Learning from local to global: An efficient distributed algorithm for
  modeling time-to-event data.
\newblock {\em Journal of the American Medical Informatics Association},
  27:1028--1036.

\bibitem[Fan et~al., 2023]{2021Communication}
Fan, J., Guo, Y., and Wang, K. (2023).
\newblock Communication-efficient accurate statistical estimation.
\newblock {\em Journal of the American Statistical Association},
  118:1000--1010.

\bibitem[Feng et~al., 2014]{2014Distributed}
Feng, J., Xu, H., and Mannor, S. (2014).
\newblock Distributed robust learning.
\newblock {\em arXiv e-prints arXiv:1409.5937}.

\bibitem[Friedman et~al., 2010]{2010Regularization}
Friedman, J.~H., Hastie, T., and Tibshirani, R. (2010).
\newblock Regularization paths for generalized linear models via coordinate
  descent.
\newblock {\em Journal of Statistical Software}, 33:1--22.

\bibitem[Huang and Huo, 2019]{huo2019distributed}
Huang, C. and Huo, X. (2019).
\newblock A distributed one-step estimator.
\newblock {\em Mathematical Programming}, 174:41--76.

\bibitem[Jordan et~al., 2019]{Jordan2019Communication}
Jordan, M.~I., Lee, J.~D., and Yang, Y. (2019).
\newblock Communication-efficient distributed statistical inference.
\newblock {\em Journal of the American Statistical Association}, 114:668--681.

\bibitem[Lamport et~al., 1982]{1982the}
Lamport, L., Shostak, R., and Pease, M. (1982).
\newblock The byzantine generals problem.
\newblock {\em ACM Transactions on Programming Languages and Systems},
  4:382--401.

\bibitem[Lecué and Lerasle, 2020]{Guillaume2020Robust}
Lecué, G. and Lerasle, M. (2020).
\newblock Robust machine learning by median-of-means: Theory and practice.
\newblock {\em The Annals of Statistics}, 48:906--931.

\bibitem[Lee et~al., 2017]{2017Communication}
Lee, J.~D., Qiang, L., Sun, Y., and Taylor, J.~E. (2017).
\newblock Communication-efficient sparse regression.
\newblock {\em The Journal of Machine Learning Research}, 18:1--30.

\bibitem[Lugosi, 2019]{2019Regularization}
Lugosi, G. (2019).
\newblock Regularization, sparse recovery, and median-of-means tournaments.
\newblock {\em Bernoulli}, 25:2075--2106.

\bibitem[Minsker, 2015]{minsker2015geometric}
Minsker, S. (2015).
\newblock Geometric median and robust estimation in banach spaces.
\newblock {\em Bernoulli}, 21:2308--2335.

\bibitem[Minsker, 2019]{minsker2019distributed}
Minsker, S. (2019).
\newblock Distributed statistical estimation and rates of convergence in normal
  approximation.
\newblock {\em Electronic Journal of Statistics}, 13:5213--5252.

\bibitem[Rosenblatt and Nadler, 2016]{Jonathan2016On}
Rosenblatt, J. and Nadler, B. (2016).
\newblock On the optimality of averaging in distributed statistical learning.
\newblock {\em Information and Inference}, 5:379--404.

\bibitem[Shang and Cheng, 2017]{2017Computational}
Shang, Z. and Cheng, G. (2017).
\newblock Computational limits of a distributed algorithm for smoothing spline.
\newblock {\em The Journal of Machine Learning Research}, 18:1--37.

\bibitem[Su and Xu, 2019]{su2019securing}
Su, L. and Xu, J. (2019).
\newblock Securing distributed gradient descent in high dimensional statistical
  learning.
\newblock {\em Proceedings of the ACM on Measurement and Analysis of Computing
  Systems}, 3:1--41.

\bibitem[Tu et~al., 2021]{tu2021variance}
Tu, J., Liu, W., Mao, X., and Chen, X. (2021).
\newblock Variance reduced median-of-means estimator for {Byzantine}-robust
  distributed inference.
\newblock {\em The Journal of Machine Learning Research}, 22:3780--3846.

\bibitem[Van~der Vaart, 2000]{2000Asymptotic}
Van~der Vaart, A. (2000).
\newblock Asymptotic statistics.
\newblock {\em Cambridge University Press}.

\bibitem[Wang et~al., 2017]{2017Efficient}
Wang, J., Kolar, M., Srebro, N., and Zhang, T. (2017).
\newblock Efficient distributed learning with sparsity.
\newblock {\em International Conference on Machine Learning}, 70:3636--3645.

\bibitem[Yin et~al., 2018]{2018Byzantine}
Yin, D., Chen, Y., Ramchandran, K., and Bartlett, P.~L. (2018).
\newblock Byzantine-robust distributed learning: Towards optimal statistical
  rates.
\newblock {\em International Conference on Machine Learning}, 80:5650--5659.

\bibitem[Yin et~al., 2019]{2019Defending}
Yin, D., Chen, Y., Ramchandran, K., and Bartlett, P.~L. (2019).
\newblock Defending against saddle point attack in byzantine-robust distributed
  learning.
\newblock {\em International Conference on Machine Learning}, 97:7074--7084.

\bibitem[Zhang et~al., 2013]{2013Communication}
Zhang, Y., Duchi, J.~C., and Wainwright, M.~J. (2013).
\newblock Communication-efficient algorithms for statistical optimization.
\newblock {\em The Journal of Machine Learning Research}, 14:3321--3363.

\bibitem[Zou and Yuan, 2008]{2008Composite}
Zou, H. and Yuan, M. (2008).
\newblock Composite quantile regression and the oracle model selection theory.
\newblock {\em The Annals of Statistics}, 36:1108--1126.

\end{thebibliography}
\end{document}